\documentclass[a4paper,12pt,oneside]{article}
\usepackage{a4wide}
\usepackage{t1enc}
\usepackage{amsfonts}
\usepackage{theorem}
\usepackage{amsmath,amssymb}
\usepackage{color}
\usepackage{graphicx}
\usepackage{psfrag}
\usepackage{ae,aecompl}
\usepackage{rotating}
\usepackage{multirow}

\newtheorem{theorem}{Theorem}[section]

\newtheorem{remark}[theorem]{Remark}

\newcommand{\hp}{\hat \varphi}
\newcommand{\Tg}{T_G}
\newcommand{\Tl}{T_L}
\newcommand{\rr}{\mathbb{R}}
\newcommand{\tor}{\mathbb{T}}
\newcommand{\T}{\mathcal{T}}
\newcommand{\wv}{\mathbf{w}}
\newcommand{\Tau}{\mathcal{T}}
\newcommand{\eps}{\varepsilon}
\newcommand{\fv}{\dot{\hat \varphi}}

\newcommand{\vect}[2]{\left ( \begin{array}{c} #1 \\ #2 \end{array} \right )}
\newcommand{\matriu}[4]{\left ( \begin{array}{cc} #1 & #2 \\ #3 & #4 \end{array} \right )}

\def\mud{\mu_{+}}
\def\mue{\mu_{-}}
\def\dosde{\mu}
\def\rhoF{\rho}

\def\ocal{\mathcal{O}}
\def\NansaFy0{\Gamma_{RD \rightarrow LD}}
\def\NansaGx0{\Gamma_{LD \rightarrow RD}}
\def\varF{h^{RL}}
\def\varG{h^{LR}}
\def\varFbar{\bar{h}^{RL}}
\def\varGbar{\bar{h}^{LR}}

\def\dotvarG{\dot{h}^{LR}}
\def\CVar{C}
\def\ustar{u^{*}}
\def\vstar{v^{*}}

\title{Quasi-periodic perturbations of heteroclinic attractor networks}

\author{Amadeu Delshams$^{1,2}$, Antoni Guillamon$^{1}$ and Gemma Huguet$^{3}$\\
\parbox{12.5cm}{
  \small
  \begin{itemize}
  \item[$^1$]
    Departament de Matem\`atiques, Universitat Polit\`ecnica de Catalunya, Avda. Dr. Mara\~{n}on 44-50, 08028 Barcelona. \\
  \item[$^2$]
    Lab of Geometry and Dynamical Systems, Universitat Polit\`ecnica de Catalunya, Avda. Dr. Mara\~{n}on 44-50, 08028 Barcelona. \\
  \item[$^3$]
    Departament de Matem\`atiques, Universitat Polit\`ecnica de Catalunya, Avda. Diagonal 647, 08028 Barcelona. \\
     \end{itemize}
}}

\date{}

\begin{document}

\maketitle

\noindent \textbf{Corresponding author:} Gemma Huguet, \texttt{gemma.huguet@upc.edu} \\



\section*{Abstract}

We consider heteroclinic attractor networks motivated by models of competition between neural populations during binocular rivalry. We show that Gamma distributions of dominance times observed experimentally in binocular rivalry and other forms of bistable perception, commonly explained by means of noise in the models, can be achieved with quasi-periodic perturbations. For this purpose, we present a methodology based on the separatrix map to model the dynamics close to heteroclinic networks with quasi-periodic perturbations. Our methodology unifies two different approaches, one based on Melnikov integrals and another one based on variational equations. We apply it to two models: first, to the Duffing equation, which comes from the perturbation of a Hamiltonian system and, second, to a heteroclinic attractor network for binocular rivalry, for which we develop a suitable method based on Melnikov integrals for non-Hamiltonian systems. In both models, the perturbed system shows chaotic behavior while dominance times achieve good agreement with Gamma distributions. Moreover, the separatrix map provides a new (discrete) model for bistable perception which, in addition, replaces the numerical integration of time-continuous models and, consequently, reduces the computational cost and avoids numerical instabilities.

\section*{Manuscript significance}
\textbf{
We study the dynamics resulting from quasi-periodic perturbations of heteroclinic attractor networks. This is a novelty with respect to previous studies which only consider periodic perturbations.
We compute explicitly a map that describes the perturbed dynamics close to the heteroclinic cycle and show that there exist chaotic orbits.
From a modelling point of view, the study of heteroclinic attractor networks was motivated by models of competition between neural populations during binocular rivalry.
We have proved that important features attributed to psychophysical experiments of bistable perception can be reproduced by quasi-periodic perturbations with two or more non-resonant frequencies.
This fact was known for noisy perturbations but not for deterministic ones.
Our methodology, based on the separatrix map, unifies two different approaches, one based on Melnikov integrals and another one on variational equations.
Moreover, we present an extension of the Melnikov approach to compute the separatrix map for non-Hamiltonian systems.
Our results provide a new (discrete) model for bistable perception, which, in addition, replaces the numerical integration of a time-continuous model.
}

\section{Introduction}

Heteroclinic networks consist of the union of several heteroclinic cycles, that is, a chain of separatrix connections between saddle points,
see \cite{Ashwin1999} for a more general and precise definition. The mathematical interest on heteroclinic cycles and networks boosted in the late 80's
(see \cite{GH1988}, \cite{ArmbrusterGH1988}, \cite{Melbourne1989}, \cite{Brannath1994}, \cite{krupa_melbourne_1995}, \cite{Krupa1997} and \cite{ArmbrusterStoneKirk2003} among others),
and soon it emerged as a suitable approach to model physical phenomena, mostly in ecology/population dynamics (see, for instance, \cite{Hofbauer1994}, \cite{Hofbauer1998}, \cite{Afraimovich2008}, \cite{AfraimovichHsuLin2008}) and more recently in neuroscience:
generation and reshaping of sequences in neural systems \cite{Rabinovich2006}, transient cognitive dynamics, metastability and decision making \cite{Rabinovich2008-PLoSCB}, decision making with memory \cite{AshwinPostlethwaite2013},
sequential memory or binding dynamics \cite{Afraimovich2015} and central pattern generators \cite{Varona2004}. Here, we focus on a specific application to cognitive neuroscience, namely, the phenomenon of bistable perception (see \cite{Ashwin2010}), as a paradigmatic example to understand the effect of quasi-periodic perturbations on heteroclinic networks.

Bistable perception consists of spontaneous alternances of sensory percepts. In humans, data is mainly obtained from psychophysical experiments that provide perceptual traces whose statistics help to understand the operating regimes of the involved brain areas. In particular, the distribution of dominance times (periods of time when a unique percept is detected, here called $\T_{dom}$) is believed to encode switching mechanisms.
As reported along the literature, see for instance \cite{FoxHerrmann67,StoneH90,Logothetis1996} and subsequent works,
such distribution of dominance times typically follows a Gamma distribution.

Models of bistable perception (see \cite{Huguet15} for a review) have two main ingredients: the presence of two states in the phase space that represent the two percepts, plus a mechanism allowing to switch from one state to another,
which can be either a negative feedback (endogenous) or noise. To account for the two states, the most well-known paradigms are oscillations \cite{Logothetis1996,Wilson2003}
and the existence of two (cross-inhibiting) point-attractors \cite{Laing2002,Moreno2007,Shpiro2009}. Ashwin and Lavric in \cite{Ashwin2010} proposed a heteroclinic attractor network model for binocular rivalry (HBR, from now on)
as another paradigm to account for the switching mechanism, inspired on previous models of winnerless competition for neural processes, see for instance \cite{Rabinovich2006,Rabinovich2008-Science,Rabinovich2008-PLoSCB}.
In these heteroclinic networks the dominance time corresponds to the time spent nearby each saddle of the heteroclinic cycle.
On the other hand, in all models noise is believed to play a more dominant role than negative feedback in shaping the trajectories of the models to fit the statistical distribution of dominance times (see \cite{Moreno2007,Pastukhov2013}).
In particular, in \cite{Ashwin2010}, the authors show that noisy perturbations lead to a Gamma-like distribution of dominance times $\T_{dom}$ for the heteroclinic attractor network.

In this paper, we use heteroclinic attractor network models to focus on the nature
of the perturbations. As explained above, models in the literature whose output fits to Gamma distributions employ noisy perturbations.
To delve into this issue, we decided to explore a minimal perturbation scheme leading to Gamma distributions of dominance times; in particular, we consider the effects of
quasi-periodic perturbations with a finite number of (non-resonant) frequencies and use appropriate tools from dynamical systems theory to study them. Besides being of mathematical interest, this issue has also a modeling relevance since it impinges onto the question of how many inputs are necessary for a specific brain area to make a perceptual decision (see also the discussion in Section \ref{sect:Discussion}).

We study heteroclinic attractor network models with periodic and, as a novelty, quasi-periodic perturbations with up to three frequencies.
Our goal is to describe the dynamics of the heteroclinic attractor network around the heteroclinic connections by means of a composition of maps between specific Poincar\'{e} sections:
a concatenation of local maps close to the saddle points together with global maps that describe the dynamics close to the heteroclinic connections. This concatenated map,
known as the \emph{separatrix map} \cite{Chirikov1979,FSZ1967,Shilnikov1965,ZF1968,PiftankinT07,AfraimovichH03}, is presented here as an alternative discrete model
for bistable perception. For this map, we carry out a thorough study of the dynamics around the heteroclinic cycle,
explore the existence of chaos by means of the computation of Lyapunov exponents and monitor the dominance (passage) times of the chaotic trajectories. Of course, studying a map avoids numerically unstable computations and increases simulation speed.

We present a methodology that (i) introduces the separatrix map as a model close to heteroclinic networks,
(ii) incorporates quasi-periodic perturbations, and (iii) unifies two different approaches, one based on Melnikov integrals and
another one based on variational equations.

We first apply our methodology to a classical model, the Duffing equation \cite{GuckenheimerH90}, which has a single saddle and a double homoclinic loop.
We use this model as a benchmark since it can be considered as the simplest version of a heteroclinic network. Once
the mathematical methodology has been settled down in the benchmarking model, we extend it to the HBR model introduced in
\cite{Ashwin2010}. For the Duffing equation, the unperturbed problem is Hamiltonian. Thus, we provide the separatrix map analytically
by means of Melnikov theory in appropriate action-angle variables, and numerically from the variational equations in the original variables.
For the HBR model, the unperturbed problem does not come from a Hamiltonian system, but
we take advantage of the exact knowledge of the heteroclinic connections to define a substitute of the action variable,
and develop a suitable method of Melnikov integrals for non-Hamiltonian systems. Again, as in the Duffing equation, we also
compute the separatrix map by means of variational equations and compare both approaches, which show good agreement.

In both cases, we develop the separatrix map for quasi-periodic perturbations and
perform a numerical study of the dynamics around the perturbed separatrices
for specific choices of the set of frequencies. For these quasi-periodic perturbations,
we show that the system exhibits chaotic behaviour around the separatrices and, remarkably, Gamma distributions of dominance times, thus showing that noise is not essential to explain such statistics.

The contents of the manuscript is organized as follows.
In Section \ref{Sect:Duffing}, we present the main concepts and tools by means of the Duffing equation. The Poincar\'{e} sections for this model are defined
in Section \ref{sect:SeparatrixMap}. In Section \ref{sect:LocalMap} we construct the
local maps nearby the saddles, whereas global maps are computed in Section \ref{sect:GlobalMap}.
By composing local and global maps, in Section \ref{sect:DuffingSeparatrix}
we then obtain two equivalent separatrix maps as a first-order approximation of the dynamics of the system around the perturbed separatrices.
In Section \ref{sec:numerical} we perform a numerical study of the dynamics of the separatrix map for a specific choice of the parameters.
In Section \ref{Sect:Heteroclinic} we carry out the same study for the HBR model introduced in
\cite{Ashwin2010}. In Section \ref{sect:AshwinVar}, we compute the separatrix map using variational equations and
in Section \ref{sect:AshwinMel} we develop the adapted Melnikov theory. Section \ref{Sect:AshwinNumeric} is devoted to numerical study.
We include three appendices with some technical details: the development of the variational equations for the Duffing system in Appendix~\ref{ap:variational},
the optimal choice of Poincar\'{e} sections in Appendix~\ref{ap:cr} and values of histogram fittings in Appendix~\ref{ap:fitting}.

\section{Dissipative separatrix map for the Duffing equation}\label{Sect:Duffing}

We use the Duffing equation as a benchmark to obtain the separatrix map both using Melnikov integrals and variational equations. Both methods use the linear local approximation, but they differ in the way to compute a \emph{global} map that will be explained in detail for each case. The Duffing equation is obtained, see \cite{GuckenheimerH90}, from the Hamiltonian
\begin{equation}\label{eq:Ham_duff}
H(x,y)=\frac{y^2}{2}-\frac{x^2}{2}+\frac{x^4}{4},
\end{equation}
where $(x,y) \in \rr^2$ and can be written as the system
\begin{equation}\label{eq:HamiltonianDuffing}
\left\{
\begin{array}{rcl}
\dot{x} & = & y, \\
\dot{y} & = & x-x^3.\\
\end{array}
\right.
\end{equation}

\begin{figure}[h]
\begin{center}
\begin{tabular}{lclcl}
(a)  & & (b) & & (c)\\
\includegraphics[width=5truecm]{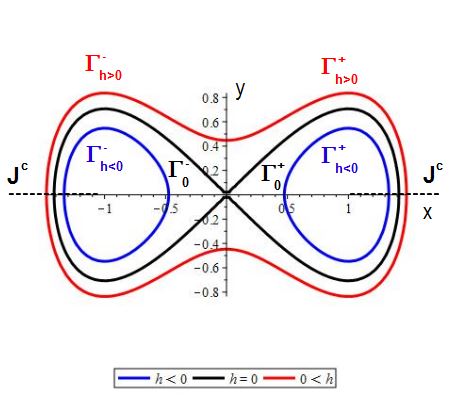} &  &
\includegraphics[width=5truecm]{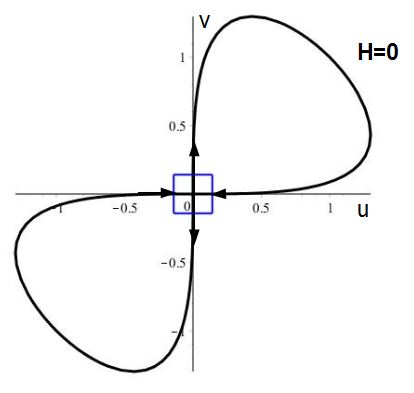} &  &
\includegraphics[width=5truecm]{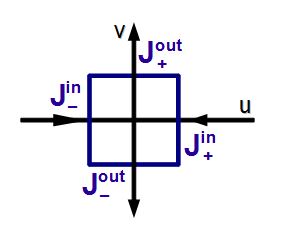}
\end{tabular}
\caption{(a) Level sets of the Hamiltonian of the Duffing equation \eqref{eq:Ham_duff} and section $J^c$ defined in \eqref{eq:DefGamma}. (b) The zero level set in the space $(u,v)$.
(c) Projections, $J_{\sigma}^{in,out}$, of the Poincar\'{e} sections defined in \eqref{eq:section_def} on the phase space $(u,v)$.
}
\label{fig:LevelCurvesDuffing}
\end{center}
\end{figure}

The flow of system \eqref{eq:HamiltonianDuffing} is organized around the double-loop separatrix $\Gamma_0$ given by $H(x,y)=0$, see Figure~\ref{fig:LevelCurvesDuffing}(a).
The level sets $\{H(x,y)=h\}$ of every energy value $h\in\rr$ can be parameterized as
\begin{equation}\label{eq:DefGamma}
\Gamma_h=\Gamma_h^{-}\cup\Gamma_h^{+}=\bigcup\limits_{\sigma=\pm}\{(x_{\sigma}(s,h),y_{\sigma}(s,h)), s\in [-T_h,T_h]\},
\end{equation}
with
\[(x_{\sigma}(0,h),y_{\sigma}(0,h))\in J^{c}:=\{(-\infty,-1)\cup(1,+\infty)\}\times\{0\},\]
$T_h>0$ and $T_0=\infty$.
For $h>0$, $\Gamma_h$ is a connected curve surrounding $\Gamma_0$, while for $h<0$, $\Gamma_h^{-}$ and
$\Gamma_h^{+}$ are two disjoint closed curves lying in different connected components of the bounded region defined by $\Gamma_0$, see again Figure~\ref{fig:LevelCurvesDuffing}(a). For $h=0$, the double-loop $\Gamma_0$ can be explicitly parameterized as

\begin{equation}\label{eq:paramh0}
x_{\sigma}(s,0)  =  \sigma \frac{\sqrt{2}}{\cosh(s)}, \qquad y_{\sigma}(s,0) = \dot{x}_{\sigma}(s,0)  =  -\sigma \frac{\sqrt{2} \sinh(s)}{\cosh^2(s)},\qquad \sigma=\pm.
\end{equation}

We will focus on a perturbed version of the Duffing equation:
\begin{equation}\label{eq:PerturbedDuffing}
\left\{
\begin{array}{rcl}
\dot{x} & = & y, \\
\dot{y} & = & x-x^3 - \gamma y + \beta x^2 y + \eps \sum_{i=1}^{n} a_i \cos(\theta_i), \\
\dot{\theta} & = & \omega, \\
\end{array}
\right.
\end{equation}
where $(x,y) \in \rr^2$, $\theta=(\theta_1,\dots,\theta_n) \in \tor^n$,
$\omega=(\omega_1,\ldots,\omega_n) \in \mathbb{R}^n$, $a_i \in \mathbb{R}$, for $i=1,\ldots,n$, and $\epsilon, \beta, \gamma \ge 0$ are small parameters.

Observe that the perturbation has two parts, one autonomous (depending on $x$ and $y$) and another non-autonomous (depending on $\theta$), which will play different roles in our study.
In fact, the autonomous part has been extensively treated in the literature, see for instance \cite{GuckenheimerH90}.
In particular, for $\eps=0$, the point $(x,y)=(0,0)$ is a saddle point with eigenvalues $\lambda_{\pm} = (-\gamma \pm \sqrt{\gamma^2+4})/2$ and the bifurcation
diagram in the $(\gamma,\beta)$ parameter space has a curve of (dissipative) homoclinic connections of the form $\beta=5/4 \,\gamma+\mathcal{O}(\gamma^2)$, see
\cite[Sect. 7.3]{GuckenheimerH90} and Remark~\ref{rem:hom_exist}.

\subsection{Poincar\'{e} sections for the separatrix map}\label{sect:SeparatrixMap}

The separatrix map was first introduced in \cite{Chirikov1979,FSZ1967,Shilnikov1965,ZF1968} as a powerful method of analysis to describe the dynamics close to a homoclinic/heteroclinic loop.
It is a singular Poincar\'{e} map (see \cite[Ch. 4]{TZ2009book}) defined on a Poincar\'{e} section near the saddle points.

In this section we describe how to construct the separatrix map for the perturbed Duffing equation~\eqref{eq:PerturbedDuffing}.
The separatrix map consists of the composition of two maps: the local map that describes the dynamics in a neighbourhood of the saddle points and the global map that
describes the dynamics in the vicinity of the homoclinic loop.

We first introduce a new coordinate system $(u,v)$ such that the linearized system around the saddle point $(0,0)$ becomes diagonal. Let
$\lambda_{-}$ and $\lambda_{+}$ be the negative and positive eigenvalues of the saddle point, respectively,
and let $v_{-}=1/\sqrt{1+\lambda_{-}^2} \,(1, \lambda_{-}) = \mue\,(-\lambda_{+},1)$ and $v_{+}=
1/\sqrt{1+\lambda_{+}^2} \, (1, \lambda_{+})=\mud\,(-\lambda_{-},1)$ be a pair of corresponding eigenvectors,
where $\mu_{+}:=\lambda_+/\sqrt{1+\lambda_+^2}$ and $\mu_{-}:=\lambda_{-}/\sqrt{1+\lambda_-^2}$ (notice that
we have used that $\lambda_{+} \lambda_{-}=-1$).
Then, the coordinate change
\begin{equation}\label{eq:coordChange}
\vect{x}{y}=\CVar \vect{u}{v},\ \mbox{where }  \CVar=\matriu{-\mue\, \lambda_{+}}{-\mud \lambda_{-}}{\mue}{\mud},
\end{equation}
transforms system \eqref{eq:PerturbedDuffing} into
\begin{equation}\label{eq:systemuv}
\left\{
\begin{array}{rcl}
\dot{u} & = & F^u(u,v,\theta,\eps)= \lambda_{-} u + f^{u}(u,v) +\eps \dfrac{\lambda_{-}}{\mue\,(\lambda_{-}-\lambda_{+})}\, \sum_{i=1}^{n} a_i \cos(\theta_i),\\
\dot{v} & = & F^v(u,v,\theta,\eps)= \lambda_{+} v + f^{v}(u,v) -\eps \dfrac{\lambda_{+}}{\mud\,(\lambda_{-}-\lambda_{+})}\, \sum_{i=1}^{n} a_i \cos(\theta_i), \\
\dot{\theta} & = & \omega, \\
\end{array}
\right.
\end{equation}
where $(f^{u}(u,v),f^{v}(u,v))$ is the transformation of the perturbation $(0,-x^3 + \beta x^2 y)$ under the change of variables \eqref{eq:coordChange} and so, it consists of a pair of homogeneous polynomials of degree $3$.

We then consider four segments $J_{\sigma}^{in}$ and $J_{\sigma}^{out}$, $\sigma\in\{+,-\}$, in a neighbourhood of the saddle point, which are transversal to the unperturbed separatrix
$\Gamma_0$ and located close to the saddle point (see Figure~\ref{fig:LevelCurvesDuffing}(b) and (c)). If we consider the
angular variable $\theta \in \tor^n$, we denote by $\Sigma_{\sigma}^{in,out}=J_{\sigma}^{in,out} \times \tor^n$ the corresponding sections in the extended phase space. Mathematically,
\begin{equation}\label{eq:section_def}
\begin{array}{rcl}
\Sigma_{\sigma}^{in} & =\{ (u=\sigma u^{*},v,\theta)\} & = \{ (\sigma,s=s_{0}^{*}(h),h,\theta)\}, \\
\Sigma_{\sigma}^{out} & =\{ (u,v=\sigma v^{*},\theta)\} & = \{ (\sigma,s=-s_{1}^{*}(h),h,\theta)\}.
\end{array}
\end{equation}
The relationship between $u^{*}, v^{*}$ and $s_{0}^{*}(h), s_{1}^{*}(h)$ is given by
\[ \vect{u^{*}}{-} = \CVar^{-1} \vect{x_{+}(s_0^{*}(h),h)}{y_{+}(s_0^{*}(h),h)} \quad
\mbox{ and } \quad
\vect{-}{v^{*}}= \CVar^{-1} \vect{x_{+}(-s_1^{*}(h),h)}{y_{+}(-s_1^{*}(h),h)}.\]

\begin{remark}
For the sake of clarity we will omit the sign $\sigma$ in the derivation of the separatrix map.
\end{remark}

\subsection{The local map}\label{sect:LocalMap}
The local map describes the dynamics from section $\Sigma^{in}:=\Sigma^{in}_{+}\cup \Sigma^{in}_{-}$ to section $\Sigma^{out}:=\Sigma^{out}_{+}\cup \Sigma^{out}_{-}$ by approximating it by the linearized dynamics around the saddle point $(0,0)$
of \eqref{eq:systemuv} for $\eps=0$. Thus,
\[\vstar=e^{\lambda_{+} \T} v + \ocal(\eps,\gamma,|u|^3+|v|^3), \qquad u=e^{\lambda_{-} \T}\ustar + \ocal(\eps,\gamma,|u|^3+|v|^3) \qquad \textrm{ and } \T=\frac{1}{\lambda_{+}} \ln \left | \frac{\vstar}{v} \right |.\]
Assuming that $|u|,|v|,\eps$ and $\gamma$ are small, we can neglect the terms $\ocal(\eps,\gamma,|u|^3+|v|^3)$ and write the local map in the variables $(u,v,\theta)$ as:

\begin{equation*}
\begin{array}{rccl}
\Tl: & \Sigma^{in} & \rightarrow & \Sigma^{out} \\
& (v,\theta) & \mapsto & (\bar{u},\bar{\theta}) \\
\end{array}
\end{equation*}
with
\begin{equation}\label{eq:linearuv}
\begin{array}{rcl}
\dfrac{\bar{u}}{u^{*}} & = & \left |\dfrac{v}{\vstar} \right |^{\nu},  \\
\bar{\theta} & = & \theta + \dfrac{\omega}{\lambda_{+}} \ln \left | \dfrac{\vstar}{v} \right |, \\
\end{array}
\end{equation}
where $\nu=-\lambda_{-}/\lambda_{+}>0$.

Equivalently, using the global variable $h= 2\,\mud \,\mue\, u\,v + \gamma \ocal_2(u,v)$,
the local map in the variables $(h,\theta)$ is given up to first order in $(u,v)$ by
$$
\begin{array}{rccl}
\Tl: & \Sigma^{in} & \rightarrow & \Sigma^{out} \\
& (h,\theta) & \mapsto & (\bar{h},\bar{\theta}) \\
\end{array}
$$
where
\begin{equation}\label{eq:DuffingLocalMap}
\begin{array}{rcl}
\dfrac{\bar{h}}{2\, \mud\, \mue\, \ustar \vstar} & = & \left (\dfrac{h}{2\, \mud\, \mue\, \ustar \vstar} \right)^{\nu}  \\
& & \\
\bar{\theta} & = & \theta + \dfrac{\omega}{\lambda_{+}} \ln \left|\dfrac{2\, \mud\, \mue\, \ustar \vstar}{h}\right|, \\
\end{array}
\end{equation}
where we have used that $h$ and $\bar{h}$ are given in first order in $(u,v)$ by $h = 2\, \mud\, \mue\, \ustar v$
and $\bar{h} = 2\, \mud\, \mue\, \bar{u} \vstar$.

\subsection{The global map}\label{sect:GlobalMap}
The global map describes the dynamics from section $\Sigma^{out}$ to section $\Sigma^{in}$ by means
of the linearized dynamics around the separatrix. In this section we will discuss two
different approaches to compute it. In one case, we will use variational equations to describe
the dynamics around the separatrix for $\beta=5/4 \, \gamma + \ocal(\gamma^2)$ and $\gamma \geq 0$ and
considering $\eps$ as the perturbation parameter. In the other case, we will use Melnikov integrals to
describe the dynamics around the unperturbed separatrix assuming that both $\eps$ and $\gamma$ (and $\beta$)
are the perturbation parameters.

\subsubsection{The global map via variational equations}\label{sect:VariationalDuffing}

In this section we will compute the global map in the variables $(u,v,\theta)$ using variational
equations. In these variables the global map is defined as
\begin{equation}\label{eq:Tgf}
\begin{array}{rccl}
\Tg: & \Sigma^{out} & \rightarrow & \Sigma^{in} \\
& (u,\theta) & \mapsto & (\bar{v},\bar{\theta}) \\
\end{array}
\end{equation}
where $(\ustar,\bar{v},\bar{\theta})=\varphi (\tau^{*}(u,\vstar,\theta); u,\vstar,\theta)$ and $\varphi(t;\wv_0)$ is the solution of system~\eqref{eq:systemuv} with initial condition $\wv_0$.

In order to obtain an approximation of the global map, we consider the local dynamics around the
separatrix $\Gamma$ that exists for $\beta=5/4 \,\gamma + \ocal (\gamma^2)$ and $\eps=0$.
Thus, we now consider the parameter
$\eps$ as a variable (i.e. $\dot{\eps}=0$) and denote by $\hat{\varphi} (t;u,v,\theta,\eps)$ the flow
of the extended system \eqref{eq:systemuv} adding $\dot{\eps}=0$, and $\hat \Tg$ the extended global map.
Moreover, let us consider the point $\wv^s:=(u^s,v^{*},\theta^s,\eps=0)$ where $u^s$ is the
$u$-coordinate of the point $p=\{ J^{out} \cap \Gamma\}$, and $\theta^s \in \mathbb{T}^n$. The
image of this point under the extended global map is the point $(u^{*},v^s,\theta^s+\omega \T^{*},0)$ where $v^s$ is the
$v$-coordinate of the point $q=\{ J^{in} \cap \Gamma\}$ and
\begin{equation}\label{eq:tau_star}
\T^{*}=\tau^{*}(\wv^s), \, \mbox{ where } \wv^s=(u^s,v^{*},\theta^s,0),
\end{equation}
which is independent of $\theta^s$ when $\eps=0$. Thus, for any point of the form
$(u,v^{*},\theta,\eps)=(u^s + \Delta u, v^{*}, \theta^s + \Delta \theta, \eps)$ its image for $\hat \Tg$ is given by
\[
\begin{array}{rcl}
\hat \Tg (u,v^{*},\theta,\eps) &=& \hat \Tg (\wv^s)  + D \hat \Tg (\wv^s) \cdot \Delta + \ocal (\Delta^2),\\
& = & (u^{*},v^s,\theta^s + \omega \T^{*},0)+ D \hat \Tg (\wv^s) \cdot \Delta + \ocal (\Delta^2),
\end{array}
\]
where $\Delta=(\Delta u, 0, \Delta \theta, \eps)$. Notice that
\begin{eqnarray*}
D \hat \Tg (\wv^s) & = & D \hat \varphi (\tau^{*}(\wv^s); \wv^s) \\
& = & D_{\wv} \hat \varphi (\tau^{*}(\wv^s); \wv^s)+ \frac{\partial \hat \varphi}{\partial t}(\tau^{*}(\wv^s); \wv^s) D_{\wv} \tau^{*}(\wv^s), \\
\end{eqnarray*}
where $\wv=(u,v,\theta,\eps)$. Of course, $D_{\wv} \hat \varphi$ can be computed by means
of solving variational equations and $D_{\wv} \tau^{*}$ can be obtained from
\[ \hat \varphi^u (\tau^{*}(\wv);\wv)=u^{*},\]
where $\hat \varphi^u$ denotes the $u$-coordinate of the extended flow $\hat \varphi$. Thus, we have
\[D_{\wv} \tau^{*} =- \left ( \frac{\partial \hat \varphi^u}{\partial t} \right)^{-1}  \, D_{\wv} \hat \varphi^u.\]
Working out the details (see Appendix~\ref{ap:variational}) we obtain that
\begin{equation}\label{eq:Tgbar}
\hat \Tg (u,v^{*},\theta,\eps)= (u^{*},v^s,\theta^s+\omega \T^{*},0)+ (0, \alpha_v \Delta u + \eps \rhoF_v(\theta),
\alpha_{\theta} \Delta u + \Delta \theta + \eps \rhoF_{\theta}(\theta),\eps)+ \ocal (\Delta^2),
\end{equation}
with
\begin{equation}\label{eq:coef_uv}
\begin{array}{cc}
\alpha_v = \hat \varphi_u^v- \dfrac{F^v}{F^u} \hat \varphi_u^u, & \rhoF_v(\theta)= \hat \varphi_{\eps}^v- \dfrac{F^v}{F^u} \hat \varphi_{\eps}^u, \vspace{0.2cm}\\
\alpha_{\theta} = - \omega \dfrac{\hat \varphi_{u}^u}{F^u}, & \rhoF_{\theta}(\theta)= - \omega\dfrac{\hat \varphi_{\eps}^u}{F^u}, \\
\end{array}
\end{equation}
where $F^u$ and $F^v$ are given in \eqref{eq:systemuv}; the subindex in $\hat \varphi$ denotes derivation with respect to that variable and the superindex denotes the corresponding coordinate. Moreover, they are obtained by means of solving
variational equations (see Appendix \ref{ap:variational}). Using that $\Delta \theta = \theta -\theta^s$ and $\Delta u = u -u^s$, and
disregarding the terms of $\ocal(\Delta^2)$ in \eqref{eq:Tgbar}, the global map $\Tg$ describing the dynamics
from $\Sigma^{out}$ to $\Sigma^{in}$
is given by
\begin{equation*}
\begin{array}{rccl}
\Tg: & \Sigma^{out} & \rightarrow & \Sigma^{in} \\
& (u,\theta) & \mapsto & (\bar{v},\bar{\theta}) \\
\end{array}
\end{equation*}
where
\begin{equation}\label{eq:exp_global}
\begin{array}{rcl}
\bar{v}  & = & v^s + \alpha_v (u-u^s) + \eps \rhoF_v(\theta), \\
\bar{\theta} & = & \theta + \omega \T^{*} + \alpha_{\theta} (u-u^s) + \eps \rhoF_{\theta}(\theta).
\end{array}
\end{equation}
Since $u-u^{s}$ and $\eps$ are assumed to be small, the contribution of the terms $\alpha_{\theta} (u-u^s)$ and $\eps \rhoF_{\theta}$ is negligible compared to the finite term $\omega \T^{*}$. Moreover, one can see that the terms $u^s$ and $v^s$ are $\ocal_2(\ustar,\vstar)$. Therefore, considering only the dominant terms,
we can write the global map as
\begin{equation}\label{eq:global_variacionals}
\begin{array}{rcl}
\bar{v} & = & \alpha u  + \eps \rhoF(\theta),\\
\bar{\theta} & = & \theta + \omega \T^{*},
\end{array}
\end{equation}
where $\alpha=\alpha_v$ and $\rho=\rho_v (\theta)$ are defined in \eqref{eq:coef_uv} and
$\Tau^{*}$ is defined in \eqref{eq:tau_star}. Notice that we have removed the subscript $v$ from $\alpha$ and $\rhoF$.

\begin{remark}
In the literature, the parameter $\alpha$ is taken to be $\alpha=1$ (see \cite{AfraimovichH03,GonchenkoSV13}). In this paper, we do not assume it to be 1 but we will compute it explicitly for some examples (see
Section~\ref{sec:numerical}).
\end{remark}

\subsubsection{The global map via Melnikov integrals}\label{sect:MelnikovDuffing}
In this section we will compute the global map using the $(h,s)$ variables by means of the
Melnikov integral. Taking advantage of the fact that the unperturbed
system \eqref{eq:HamiltonianDuffing} is Hamiltonian, we have that
\[
\begin{array}{rcl}
\dot{h} &= & H_y\, q, \\
\dot{s} &= & 1+s_y\, q, \\
\dot{\theta} &= & \omega, \\
\end{array}
\]
where $H$ is the Hamiltonian function in \eqref{eq:Ham_duff} and
$q(x,y,\theta)=-\gamma y + \beta x^2 y + \eps \sum_{i=1}^{n} a_i \cos(\theta_i)$
is the perturbation in \eqref{eq:PerturbedDuffing}. All the functions of the above expression are evaluated on $x=x(s,h)$ and $y=y(s,h)$, introduced in \eqref{eq:DefGamma}.

Therefore, the global map describing the dynamics from $\Sigma^{out}$ to $\Sigma^{in}$
is defined as:
\[
\begin{array}{rccl}
\Tg: & \Sigma^{out} & \rightarrow & \Sigma^{in} \\
& (h,\theta) & \mapsto & (\bar{h},\bar{\theta}) \\
\end{array}
\]
where, using that $s_i^{*}(h)=s_i^{*}(0)+\ocal(h)$, for $i=0,1$, we have
\begin{equation}\label{eq:DuffingGlobalMap}
\begin{array}{rcl}
\bar{h} - h & = & \int_{-s^{*}_{1}}^{s^{*}_{0}} H_y\, q = M(\theta) + (\gamma+\beta+\eps)\,\ocal(u^{*},v^{*},|h|,\eps,\gamma,\beta) \\
\bar{\theta} - \theta & = & \omega (s^{*}_{0} + s^{*}_{1}) + \ocal(|h|,\gamma,\eps,\beta), \\
\end{array}
\end{equation}
where $s^{*}_i$ denotes $s^{*}_i(0)$ and $M(\theta)$ is the Melnikov integral
for system \eqref{eq:PerturbedDuffing} on the level curve $H(x,y)=h=0$.

We compute the Melnikov integral on the positive branch of the level curve $H(x,y)=0$, i.e.
$M(\theta)=M_{+}(\theta)$ (the case for the negative branch is
analogous). Let us denote $x_0(t)=x_{+}(t,0)$ and $y_0(t)=\dot{x}_0(t)=y_{+}(t,0)$
the parameterization of $\Gamma_0^{+}$ given in
equation~\eqref{eq:paramh0}. Then,
\begin{equation*}
\begin{array}{rcl}
M(\theta) &= & \displaystyle \int_{-\infty}^{\infty} \frac{\partial H}{\partial y} (x_0(t),y_0(t))
\, q(x_0(t),y_0(t),\theta+\omega t) dt \\
& = & \displaystyle\int_{-\infty}^{\infty} \dot{x}_0 (t) \left (- \gamma \dot{x}_0 (t) + \beta x_0^2(t) \dot{x}_0 (t) + \eps \sum_{i=1}^{n} a_i \cos (\theta_i + \omega_i t)  \right) dt \\ & & \\
& = & - \gamma \displaystyle\int_{-\infty}^{\infty} \dot{x}_0^2(t)\, dt + \beta \displaystyle\int_{-\infty}^{\infty} x_0^2(t) \dot{x}_0^2(t) \, dt- \eps \sum_{i=1}^{n} a_i \sin (\theta_i) \displaystyle \int_{-\infty}^{\infty} \dot{x}_0(t) \sin (\omega_i t) dt,
\end{array}
\end{equation*}
where in the last term we have used that the integral on $(-\infty,\infty)$ of an odd function is zero.

We compute each integral separately. Thus,
\[
\int_{-\infty}^{\infty} \dot{x}_0^2(t)dt = 2 \int_{-\infty}^{\infty} \frac{\sinh^2(t)}{\cosh^4(t)}dt = \frac{4}{3}, \quad
\int_{-\infty}^{\infty} x_0^2(t) \dot{x}_0^2(t)dt = 4 \int_{-\infty}^{\infty} \frac{\sinh^2(t)}{\cosh^6(t)}dt = \frac{16}{15},
\]
and
\begin{equation*}
\begin{array}{rl}
\displaystyle \int_{-\infty}^{\infty} \sin (\omega_i t) \dot{x}_0 (t) dt &=  \left . \sin (\omega_i t) x_0 (t) \right ]_{-\infty}^{\infty} - \omega_i  \displaystyle \int_{-\infty}^{\infty} \cos (\omega_i t) x_0 (t) dt \\
& = \displaystyle - \omega_i \sqrt{2} \int_{-\infty}^{\infty} \frac{\cos (\omega_i t)}{\cosh(t)} dt  = \displaystyle - \frac{\sqrt{2} \omega_i \pi}{\cosh (\pi \omega_i/2)},
\end{array}
\end{equation*}
where the last integral has been computed by the Residue Theorem (see \cite{DelshamsG00}).
Thus, the Melnikov integral has the form
\begin{equation}\label{eq:mel_int}
M(\theta)=M_{+}(\theta)= - \frac{4}{3} \gamma + \frac{16}{15} \beta + \sqrt{2} \pi \eps \sum_{i=1}^n \frac{a_i \omega_i \sin (\theta_i)}{\cosh (\pi \omega_i /2)}.
\end{equation}

\begin{remark}\label{rem:hom_exist}
Note that $M(\theta)=0$ for $\eps=0$ and $\beta=5/4 \, \gamma$. Therefore system~\eqref{eq:PerturbedDuffing} has a
homoclinic orbit for $\eps=0$ and $\beta=5/4 \, \gamma + \ocal (\gamma^2)$, see also \cite[Sect. 7.3]{GuckenheimerH90}.
\end{remark}

\subsection{The separatrix map}\label{sect:DuffingSeparatrix}

We describe how to construct the separatrix map for the perturbed Duffing equation~\eqref{eq:PerturbedDuffing} both via variational equations and Melnikov integrals.

\subsubsection{The separatrix map via variational equations}\label{sec:ss_ssve}

We combine the local and global maps given in \eqref{eq:linearuv} and \eqref{eq:global_variacionals}, respectively,
to obtain explicit formulas for the separatrix map in the variables $(u,\theta)$.

Notice that when $v>0$, the local map takes points from $\Sigma^{in}$ to $\Sigma^{out}_{+}$, but
when $v<0$, the local map takes points from $\Sigma^{in}$ to $\Sigma^{out}_{-}$. Recovering the
variable $\sigma$ in \eqref{eq:section_def}, we can write the separatrix map as
\[
\begin{array}{rccl}
S:= \Tl \circ \Tg : & \Sigma^{out} & \rightarrow & \Sigma^{out} \\
& (u,\theta,\sigma) & \mapsto & (\bar{u},\bar{\theta},\bar{\sigma})\\
\end{array}
\]
where
\begin{equation}\label{eq:map_var}
\begin{array}{rcl}
\dfrac{\bar{u}}{u^*} & = & \sigma |v^{*}|^{-\nu} |\alpha u + \eps \rhoF(\theta)|^{\nu}, \\
\bar{\theta} & = & \theta +  \omega \T^{*} + \dfrac{\omega}{\lambda_{+}}\ln \displaystyle \left | \frac{v^{*}}{\alpha u  + \eps \rhoF(\theta)} \right | \quad \mod 2 \pi,\\
\bar{\sigma} & = & \textrm{sign}(\alpha u  + \eps \rhoF(\theta)).\\
\end{array}
\end{equation}

By scaling the variables $u$ and $v$ by $\ustar$ and $\vstar$, respectively, i.e. $u:=u/\ustar$ and $v:=v/\vstar$, the global map becomes
\begin{equation*}\label{eq:map_var_scaled}
\begin{array}{rcl}
\bar{u} & = & |\tilde{\alpha} u + \eps \tilde{\rhoF}(\theta)|^{\nu}, \\
\bar{\theta} & = & \theta +  \omega \T^{*} + \dfrac{\omega}{\lambda_{+}}\ln \displaystyle \left | \frac{1}{\tilde{\alpha} u  + \eps \tilde{\rhoF}(\theta)} \right | \quad \mod 2 \pi,\\
\bar{\sigma} & = & \textrm{sign}(\tilde{\alpha} u  + \eps \tilde{\rhoF}(\theta)),\\
\end{array}
\end{equation*}
where $\tilde{\rhoF}(\theta)=\displaystyle\frac{1}{\vstar} \rhoF(\theta)$ and $\tilde{\alpha}:=\displaystyle\frac{\ustar}{\vstar} \alpha$.
Notice that if we choose $u^*=v^*$ then $\tilde{\alpha}=\alpha$.

\subsubsection{The separatrix map via Melnikov integrals}\label{sect:SepMelnikov}

We combine the local and global maps given in \eqref{eq:DuffingLocalMap} and \eqref{eq:DuffingGlobalMap}, respectively,
to obtain explicit formulas for the separatrix map in the variables $(h,\theta)$.

Notice first that the local map takes points from $\Sigma^{in}_+$ (resp. $\Sigma^{in}_{-}$) to $\Sigma^{out}_{+}$ or $\Sigma^{out}_{-}$ depending on the value
of $h$ at the section $\Sigma^{in}_+$ (resp. $\Sigma^{in}_{-}$).
Thus, recovering the variable $\sigma$ in \eqref{eq:section_def} and using the dominant terms given in \eqref{eq:DuffingLocalMap} and \eqref{eq:DuffingGlobalMap} we can write the separatrix map as
\[
\begin{array}{rccl}
S:= \Tl \circ \Tg : & \Sigma^{out} & \rightarrow & \Sigma^{out} \\
& (h,\theta,\sigma) & \mapsto & (\bar{h},\bar{\theta},\bar{\sigma})\\
\end{array}
\]
where
\begin{equation}\label{eq:map_melnikov}
\begin{array}{rcl}
\dfrac{\bar{h}}{2 \mud\, \mue\, \ustar \vstar} & = & \left (\dfrac{h+M_{\sigma}(\theta)}{2 \mud\, \mue\, \ustar \vstar} \right )^{\nu}, \\ & & \\
\bar{\theta} & = & \theta + \omega (s_0^{*} + s_1^{*}) + \dfrac{\omega}{\lambda_{+}} \ln \left|\dfrac{2 \mud\, \mue\, \ustar \vstar}{h+M_{\sigma}(\theta)}\right|, \\
\bar{\sigma} & = & - \sigma\, \textrm{sign} (h + M_{\sigma}(\theta)).
\end{array}
\end{equation}

By scaling the variable $h$ and redefining $h:=h/(2 \mud \mue \ustar \vstar)$ the global map can be simply written as
\begin{equation*}\label{eq:map_melnikov_scaled}
\begin{array}{rcl}
\bar{h} & = & \left  (h+ \widetilde{M}_{\sigma}(\theta \right )^{\nu}, \\
\bar{\theta} & = & \theta + \omega (s_0^{*} + s_1^{*}) + \dfrac{\omega}{\lambda_{+}} \ln \left|\dfrac{1}{h+\widetilde{M}_{\sigma}(\theta)}\right|, \\
\bar{\sigma} & = & - \sigma\, \textrm{sign} (h + \widetilde{M}_{\sigma}(\theta)),
\end{array}
\end{equation*}
where $\widetilde{M}_{\sigma}(\theta):=\dfrac{1}{2 \mud\, \mue\, \ustar \vstar} M_{\sigma}(\theta)$.

\subsection{Numerical computations}\label{sec:numerical}

We compute numerically the separatrix map for the perturbed Duffing equation using variational equations and
we refer to Section~\ref{sec:differences} for the differences between the separatrix map obtained using this
approach and the one obtained via Melnikov integrals.

We consider a quasi-periodic perturbation consisting of at most 3 frequencies given by $\omega_1=1$, $\omega_2=\frac{\sqrt{5}-1}{2}$ and $\omega_3=\sqrt{769}-27$.
We choose them in order to pick three frequencies that are as much incongruent as possible.
Indeed, $\omega_3=\sqrt{769}-27$ is the real number that provides the best constant $C=C(\omega_3)\approx 0.233126\dots$ in the inequalities
\begin{equation}\label{eq:diophantine}
|k_0 + k_1\omega_2 + k_2\omega_3|\geq C (|k_0|+|k_1|+|k_2|)^{-2},
\end{equation}
for all integers satisfying $|k_1|,|k_2| \leq K:= 2^{20}$, $k_0 \in \mathbb{Z}$, amongst all numbers
$\omega_3$ that are the decimal part of numbers of the form $\sqrt{p}$, with $p$ being a prime number smaller than $1000$ \cite{Simo16}.
An alternative to these frequencies is to choose $\omega_2=\Omega$, where $\Omega$ is ``the cubic golden number'', i.e.
the real cubic root of $x^ 3 + x - 1=0$, $\Omega \approx 0.6823$ and $\omega_3= \Omega ^2$. In this case, $C=2(5+\Omega+4\Omega ^2)/31\approx 0.4867$ in the expression
\eqref{eq:diophantine} \cite{DelshamsGG14}. Simulations with these frequencies do not show significant differences
(results not shown).

We first compute the separatrix map for two orders of magnitude of $\gamma$ (and $\beta$, since $\beta$
is chosen as $\beta=5/4 \, \gamma + \ocal (\gamma^2)$ so that the map for $\eps=0$ has a homoclinic orbit, see
remark~\ref{rem:hom_exist}).
In both cases, we introduce the parameter $r$ such that $u^*=v^*=r$ and we chose $r=0.1$
(see Appendix~\ref{ap:cr} for a justification).
The coefficients of the map \eqref{eq:global_variacionals} are computed according to the formulas given in \eqref{eq:coef_uv}, which involve solving numerically the variational equations
around the separatrix (see Appendix~\ref{ap:variational}). For the numerical integration we have used a Runge-Kutta method of order
7/8 with a fixed tolerance of $10^{-12}$.

For $\gamma=0.008$, the map~\eqref{eq:map_var} is given by
\begin{equation}\label{eq:dms_g8-3}
\begin{array}{rcl}
\bar{u} &= &s r^{1-\nu} |\alpha u + \eps \rho(\theta)|^{\nu}, \\
\bar{\theta}_i & = & \theta_i + \omega_i\, 7.3752858056 + \omega_i/0.9960080000\, \ln(r/|v|), \, \textrm{for } \, i=1,2,3,\\
\bar{\sigma} & = & \textrm{sign}(v), \\
\end{array}
\end{equation}
where $\alpha = 0.9733201532$, $r=0.1$ and
\begin{equation*}
\begin{array}{rcl}
\rhoF(\theta) &= & a_1 (9.7591847996 \cos(\theta_1)+15.6872106985 \sin(\theta_1)) \\
& & +a_2 (-13.2851558002 \cos(\theta_2)+11.5920512181 \sin(\theta_2)) \\
& & +a_3 (-7.9272168789 \cos(\theta_3)+ 16.9623679354 \sin(\theta_3)), \\
\end{array}
\end{equation*}
with $\theta=(\theta_1,\theta_2,\theta_3) \in \mathbb{T}^3$ and $a_1,a_2,a_3 \in \mathbb{R}$.

For $\gamma=0.08$, the same value as in \cite{StoneH90}, the map is given by
\begin{equation}\label{eq:dms_g8-2}
\begin{array}{rcl}
\bar{u} &= &s r^{1-\nu} |\alpha u + \eps \rho (\theta)|^{\nu}, \\
\bar{\theta}_i & = & \theta_i + \omega_i\, 7.3784656185 + \omega_i/0.9607996803\, \ln(r/|v|), \, \textrm{for }\, i=1,2,3,\\
\bar{\sigma} & = & \textrm{sign}(v), \\
\end{array}
\end{equation}
where $\alpha = 0.7629736972$, $r=0.1$ and
\begin{equation*}
\begin{array}{rcl}
\rhoF(\theta) &= & a_1 (9.9901759770 \cos(\theta_1)+13.0767449862 \sin(\theta_1)) \\
& & +a_2 (-10.9333035475\cos(\theta_2)+11.2761757850\sin(\theta_2)) \\
& & +a_3 (-5.7535147048 \cos(\theta_3)+15.6518248196\sin(\theta_3)), \\
\end{array}
\end{equation*}
with $\theta=(\theta_1,\theta_2,\theta_3) \in \mathbb{T}^3$ and $a_1,a_2,a_3 \in \mathbb{R}$.
Notice that the coefficient $\alpha$ in \eqref{eq:dms_g8-3} for $\gamma=0.008$ is larger than in \eqref{eq:dms_g8-2}
for $\gamma=0.08$, showing that the
separatrix becomes more contractive as $\gamma$ increases (see Section~\ref{sec:differences}).

\begin{figure}[h]
\begin{center}
\begin{tabular}{clll}
&(a) 1 frequency & (b) 2 frequencies &  (c) 3 frequencies \\
\begin{sideways}$\qquad \gamma=0.008$\end{sideways} &
\includegraphics[width=5cm,height=3.8cm]{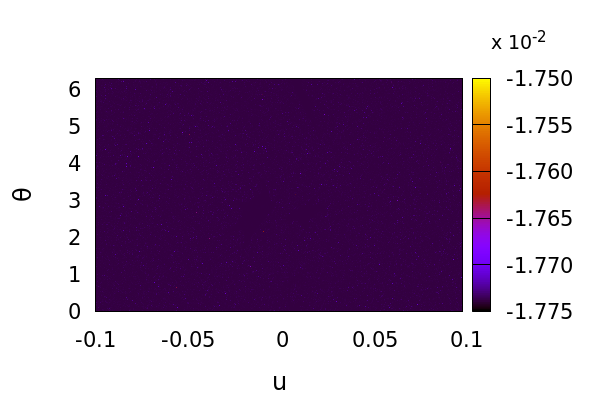} &
\includegraphics[width=5cm]{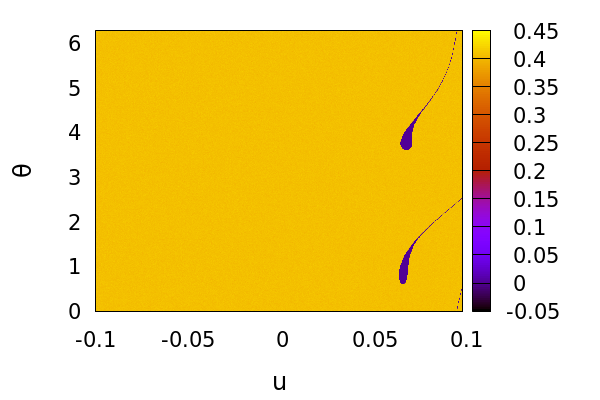} &
\includegraphics[width=5cm]{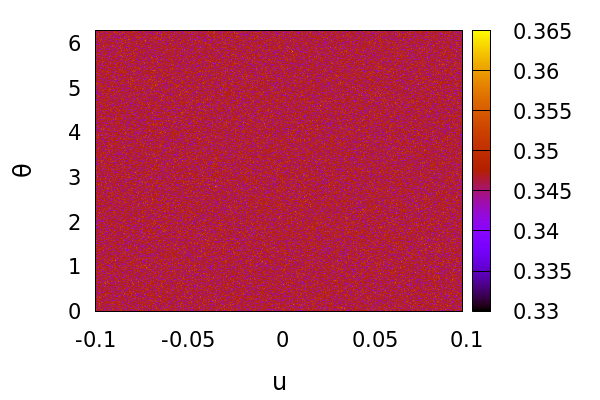} \\
\begin{sideways}$\qquad \gamma=0.08$\end{sideways} &
\includegraphics[width=5cm,height=3.8cm]{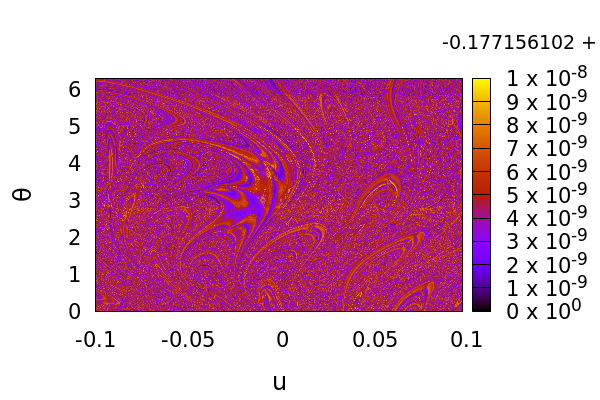} &
\includegraphics[width=5cm]{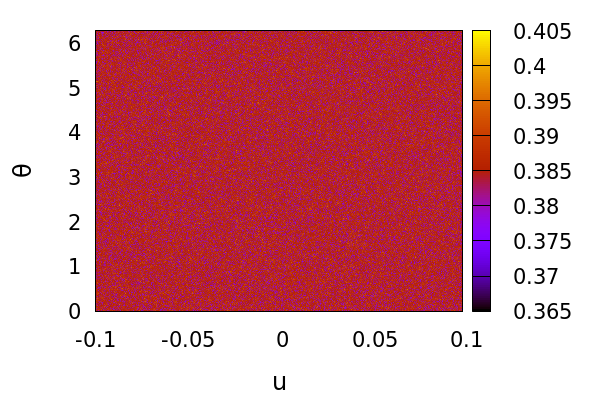} &
\includegraphics[width=5cm]{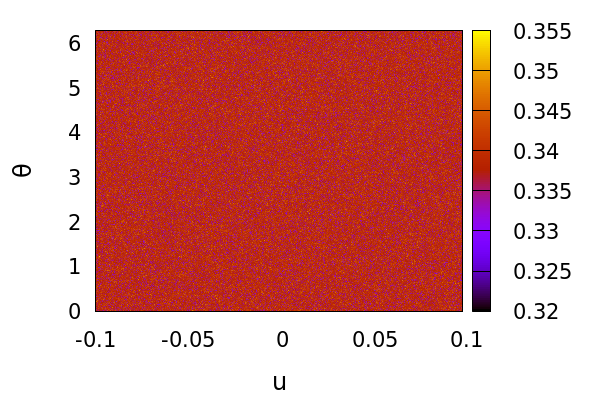} \\
\end{tabular}
\caption{Maximal Lyapunov exponent (computed using MEGNO) for the orbits of several separatrix maps of the Duffing equation with initial conditions on the phase space $(u,\theta)$ ($\theta_i=\theta$, for $i=1,2,3$):
(top) separatrix map \eqref{eq:dms_g8-3} corresponding to $\gamma=0.008$ with $r=0.1$ and $\eps=0.001$, and (bottom) separatrix map \eqref{eq:dms_g8-2} corresponding to $\gamma=0.08$
with $r=0.1$ and $\eps=0.001$. The number of frequencies in the perturbation is indicated in each panel.}
\label{fig:megno}
\end{center}
\end{figure}

\begin{figure}
\begin{center}
\begin{tabular}{ll}
(a) & (b) \\
\includegraphics[width=8cm]{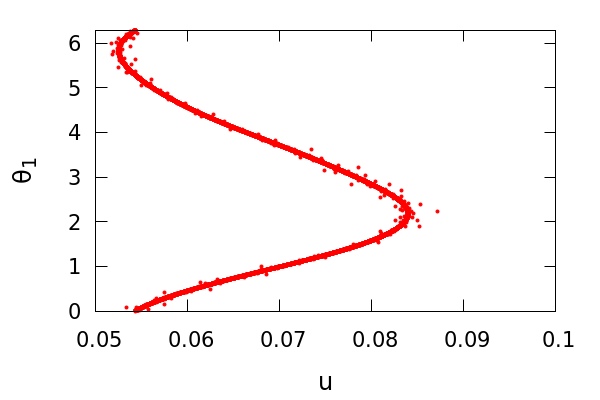} &
\includegraphics[width=8cm]{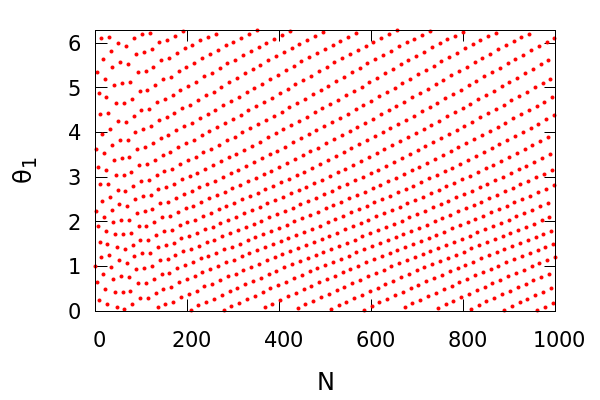} \\
\end{tabular}
\caption{Iterates of the separatrix map \eqref{eq:dms_g8-3} corresponding to $\gamma=0.008$ with $r=0.1$, $\eps=0.001$, and a perturbation with two frequencies,
with initial conditions in the non-chaotic region (see Figure~\ref{fig:megno}(b) top). (a) Iterates on the $(u,\theta_1)$ space, and (b) iterate number vs angle variable $\theta_1$.}
\label{fig:iteratsnochaotics}
\end{center}
\end{figure}

In order to understand the dynamics of these maps we have carried out a numerical exploration of the Lyapunov exponents for different orbits of the system with
different initial conditions. We run the MEGNO program \cite{CincottaG16, CincottaGS03} to compute the maximal Lyapunov exponent for the maps \eqref{eq:dms_g8-3}
and \eqref{eq:dms_g8-2} with $\eps=0.001$ and three different perturbations:  a periodic perturbation ($a_1=1, a_2=0, a_3=0$), a quasi-periodic perturbation with 2 frequencies ($a_1=1,a_2=1,a_3=0$) and
a quasi-periodic perturbation with 3 frequencies ($a_1=1,a_2=1,a_3=1$). Results for the maximal Lyapunov exponent for different initial
conditions are shown in Figure~\ref{fig:megno}.
For the case of a periodic perturbation (1 frequency), the Lyapunov exponent is negative and both maps show non-chaotic behaviour for all the initial conditions tested
(see Figure~\ref{fig:megno}(a)). Indeed the orbits of the system tend to a periodic orbit, and in the case of the map with $\gamma=0.008$,
we observed several limiting periodic orbits (results not shown) for different initial conditions. This explains why the range of the Lyapunov exponents computed for $\gamma=0.008$ is larger than the one for $\gamma=0.08$.
For the case of quasi-periodic perturbations with 2 or 3 frequencies, the Lyapunov exponent is positive for all the initial conditions tested for the map
with $\gamma=0.08$
and a large domain of the initial conditions tested for the map $\gamma=0.008$, thus showing chaotic behaviour
(see Figure~\ref{fig:megno}(b) and (c)). Notice that there is a small area of initial conditions for which the map corresponding to
$\gamma=0.008$ and a quasi-periodic perturbation with two frequencies shows non-chaotic behaviour (see
Figure~\ref{fig:megno}(b) top). The iterates for an initial condition in this non-chaotic region are shown in Figure~\ref{fig:iteratsnochaotics},
where the numerical exploration shows the existence of an attracting invariant curve.

\begin{figure}[h]
\begin{center}
\includegraphics[width=17cm]{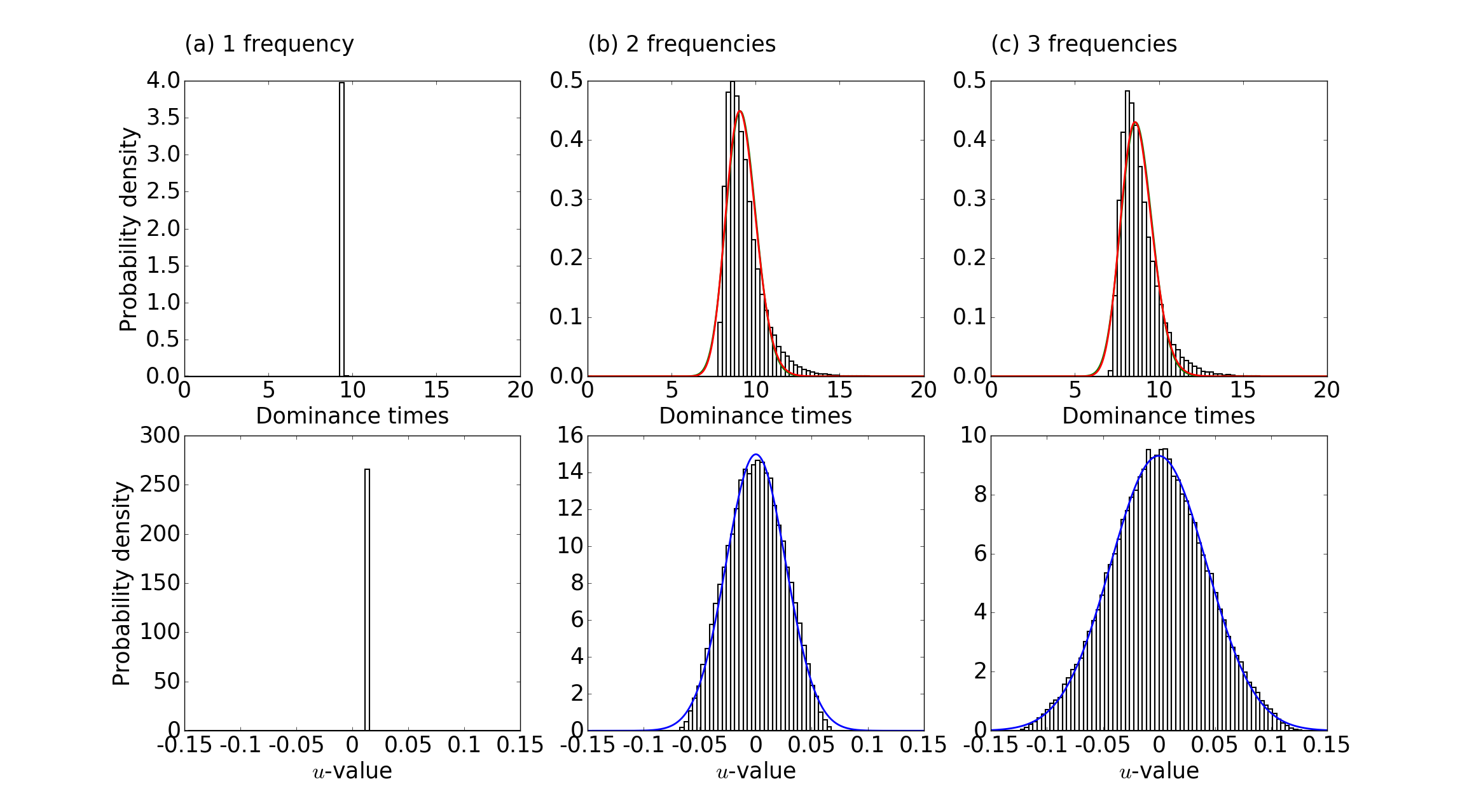}
\caption{Histograms of the dominance times (time to return to the Poincar\'{e} section $\Sigma^{out}_{v=\pm r}$) and the $u$-value at
section $v=\pm r$ for the map in \eqref{eq:dms_g8-3} corresponding to $\gamma=0.008$ with a perturbation with (a) 1, (b) 2 and (c) 3 frequencies, $\eps=0.001$, $r=0.1$ and
initial conditions $u=0$, $\theta_i=0$, $\sigma=1$ and $\gamma=0.008$.
We have used 100,000 iterates. Dominance times histograms (top) are fitted to Gamma (green) and log-normal (red) distributions (fittings not distinguishable) and
impact histograms (bottom) are fitted to normal distributions (blue), see Appendix~\ref{ap:fitting}.
Notice here that the $u$-values of the iterates for the quasi-periodic perturbation of 3 frequencies impact, in some cases, outside $r=0.1$. This is not a problem since the map computes the times correctly.
If the initial condition for the case with 1 frequency is chosen as $u=0.1$, the orbit converges to a different periodic orbit
and we obtain a histogram which is also a delta function but centered at a different position (results not shown). Histograms have been normalized to have area 1.}
\label{fig:hist_gae-3}
\end{center}
\end{figure}

\begin{figure}[h]
\begin{center}
\includegraphics[width=17cm]{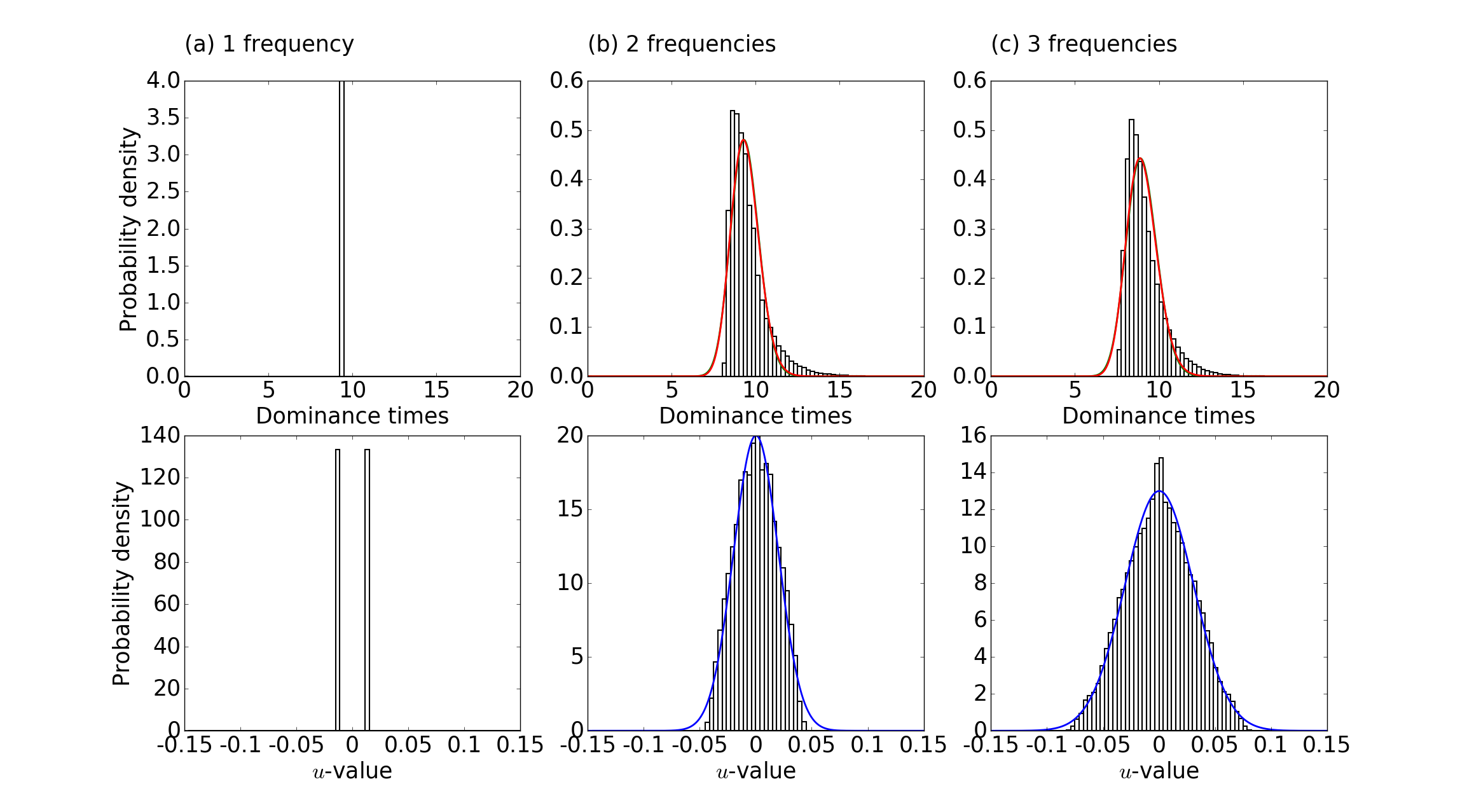}
\caption{Histograms of the dominance times (time to return to the Poincar\'{e} section $\Sigma^{out}_{v=\pm r}$) and the $u$-value at
section $v=\pm r$ for the map in \eqref{eq:dms_g8-2} corresponding to $\gamma=0.08$ with a perturbation with (a) 1, (b) 2 and (c) 3 frequencies, $\eps=0.001$, $r=0.1$ and
initial conditions $u=0$, $\theta_i=0$, $\sigma=1$ and $\gamma=0.008$.
We have used 100,000 iterates. Dominance times histograms (top) are fitted to Gamma (green) and log-normal (red) distributions (fittings not distinguishable) and
impact histograms (bottom) are fitted to normal distributions (blue), see Appendix~\ref{ap:fitting}. Histograms have been normalized to have area 1.
}
\label{fig:hist_gae-2}
\end{center}
\end{figure}

In the cases for which we detected chaotic behaviour for the orbits of the system, we explored the distribution of
the dominance times $\T_{dom}$ corresponding to the time intervals between impacts on sections
$\Sigma^{out}_{\pm}$. Mathematically, $\T_{dom}=(\bar \theta_i - \theta_i)/\omega_i$, which is independent of $i$.
In addition, we also computed the distribution of impacts on the sections $\Sigma^{out}_{\pm}$ (given by the values of $u$ along
the orbit). In particular, we considered initial conditions $u=0$, $\theta_i=0$, for $i=1,2,3$ and $\sigma=+1$ and computed the
corresponding iterates for both separatrix maps \eqref{eq:dms_g8-3} and \eqref{eq:dms_g8-2} and the three different perturbations.
In Figures~\ref{fig:hist_gae-3} and \ref{fig:hist_gae-2} we show the corresponding histograms of the dominance
times and the impacts on the section $\Sigma^{out}$.
Notice that for the cases where we observed chaotic behaviour the histograms show a log-normal or Gamma
distribution for time differences $\T_{dom}$ and a normal distribution for impacts
(see Figures~\ref{fig:hist_gae-3} and \ref{fig:hist_gae-2}, (b) and (c)),
while for the case of one frequency the histograms are just a delta function for both
distributions (see Figures~\ref{fig:hist_gae-3}(a) and \ref{fig:hist_gae-2}(a)).

\begin{figure}[h]
\begin{center}
\includegraphics[width=17cm]{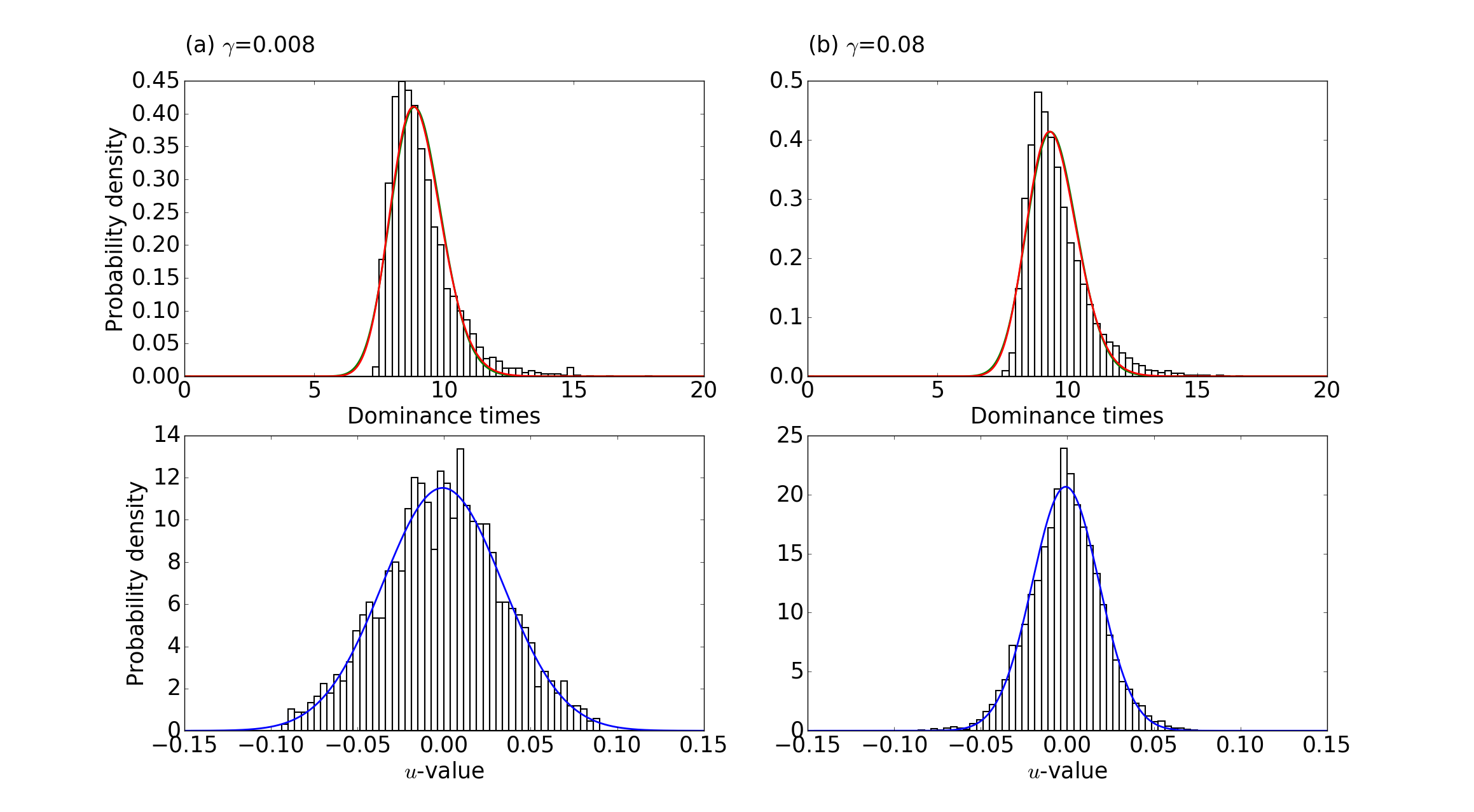}
\caption{
Histograms of the dominance times (time to return to the Poincar\'{e} section $\Sigma^{out}_{v=\pm r}$) and the $u$-value at
section $v=\pm r$ obtained from the integration of the system of differential equations
\eqref{eq:noise_uv} with $\gamma=0.008$ and $\eps=5 \cdot 10^{-4}$ (left) and $\gamma=0.08$ and $\eps=10^{-3}$ (right).
Time histograms (top) are fitted to Gamma (green) and log-normal (red) distributions (fittings not distinguishable) and
impact histograms (bottom) are fitted to normal distributions (blue), see Appendix~\ref{ap:fitting}.
We have used 10,000 iterates and initial conditions $u=0$ and $v=r$. We have integrated the system using an Euler-Maruyama method with
a stepsize of $\Delta t=10^{-6}$ for $\gamma=0.008$ and $\Delta t=10^{-5}$ for $\gamma=0.08$. For $\gamma=0.008$ the
homoclinic is less contractive, therefore the noise is chosen smaller to avoid that trajectories drift away from the separatrix. Histograms have been normalized to have area 1.}
\label{fig:hist_noise_uv}
\end{center}
\end{figure}

\begin{figure}[h]
\begin{center}
\begin{tabular}{ll}
(a) $\gamma=0.008$ & (b) $\gamma=0.08$\\
\includegraphics[width=8cm]{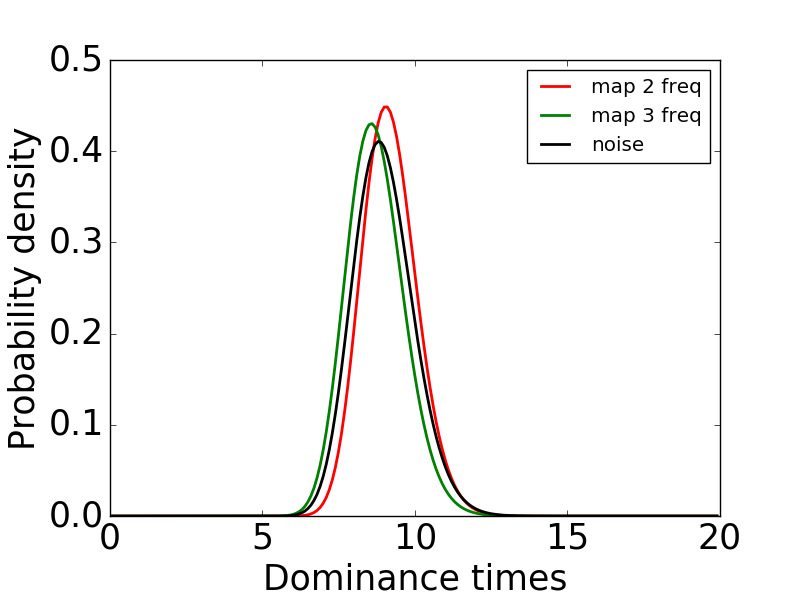} &
\includegraphics[width=8cm]{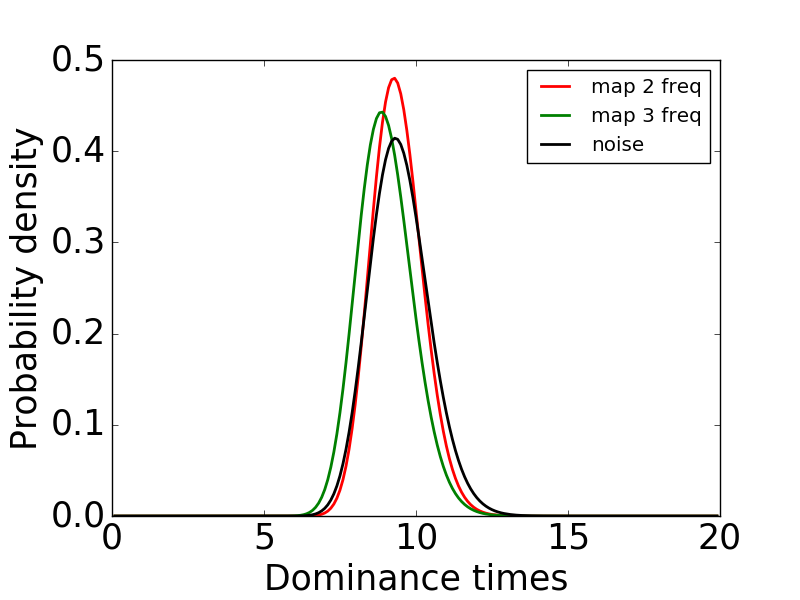} \\
\end{tabular}
\caption{Comparison between fittings to log-normal distribution of the dominance times histograms in Figures~\ref{fig:hist_gae-3}, \ref{fig:hist_gae-2} and \ref{fig:hist_noise_uv}  (similar results for fittings to Gamma distribution).}
\label{fig:resum_hist_duffing}
\end{center}
\end{figure}

We compare the results for histograms with those obtained with a perturbation consisting of white noise instead of a periodic
or quasi-periodic function.
Thus, we consider the following system of stochastic differential equations:
\begin{equation}
\begin{array}{rcl}\label{eq:noise_uv}
du & = & \lambda_{-} u + f^{u}(u,v) du  + \eps_u dW_u, \\
dv & = & \lambda_{+} v + f^{v}(u,v) dv + \eps_v dW_v,\\
\end{array}
\end{equation}
where $f^u,f^v$ are defined in \eqref{eq:systemuv} and $dW_u, dW_v$ are zero mean, independent Wiener processes. We computed the values of $u$ and $\T_{dom}$ on the section $\Sigma^{out}$ and obtained the histograms shown in Figure~\ref{fig:hist_noise_uv}. Considering a noisy perturbation in the original variables $(x,y)$ leads to similar results
(not shown), see also \cite{StoneH90}.
Notice that the shape of histograms follows a log-normal or Gamma distribution for time histograms and a normal distribution for
impacts (see Figure~\ref{fig:hist_noise_uv}). Notice also that, as observed for the case of a quasi-periodic perturbation, the case $\gamma=0.008$ contracts in a weaker way.
In order to properly compare the histograms with the ones obtained with quasi-periodic perturbations, we fit the histograms
with a log-normal and Gamma distribution (see Appendix~\ref{ap:fitting}). We show the fittings to a log-normal distribution altogether in
Figure~\ref{fig:resum_hist_duffing}. The same fittings are obtained for a Gamma distribution. We clearly see that for 2 and 3 frequencies, the histograms are similar to the ones obtained with
noise.

\subsubsection{Comparison between separatrix maps for the Duffing equation}\label{sec:differences}

Using Melnikov integrals we have obtained the separatrix map given in \eqref{eq:map_melnikov} which has an analytical expression,
while using variational equations about the perturbed separatrix ($\gamma \neq 0$, $\beta=5/4 \, \gamma + \ocal (\gamma^2)$)
we have obtained the separatrix map \eqref{eq:map_var}
which requires numerical computations. We want to compare the results obtained numerically with those obtained analytically by
means of two different methods.
We know that both maps are equivalent up to $\ocal(\gamma)$.

Let us express the map \eqref{eq:map_var} in terms of $h$. Disregarding higher order terms in $(u,v)$ we have
$h=\dosde u v^{*}$, where $\dosde=2 \mud \mue = -1 / (\lambda_{+}-\lambda_{-})$, and replacing $u$ by $h/(\mu v^{*})$ in system \eqref{eq:map_var} we have
\begin{equation}\label{eq:mapht}
\begin{array}{rl}
\dfrac{\bar h}{\mu u^{*} v^{*}}&= \left ( \dfrac{\alpha (u^{*}/v^{*}) h + \eps \dosde u^{*} \rhoF(\theta)}{\dosde u^{*} v^{*}} \right )^{\nu},\\
\bar \theta &= \theta+\omega \Tau^* + \dfrac{\omega}{\lambda_{+}} \ln \left | \dfrac{\dosde v^{*} u^{*}}{ \alpha (u^{*}/v^{*}) h + \eps \dosde u^{*} \rhoF(\theta)} \right |, \\
\bar \sigma & = - \sigma \textrm{sign} (h + \eps \mu \rho(\theta)),
\end{array}
\end{equation}
where in the last expression we have used that the sign of $v$ is the same as (resp. the opposite of) the sign of $h$ in $\Sigma^{in}_{-}$ (resp. $\Sigma^{in}_{+})$.

Now we are going to compare the map in the variables $(h,\theta)$ given in \eqref{eq:map_melnikov} with equations \eqref{eq:mapht}. First notice that to obtain
the map \eqref{eq:map_melnikov} using the Melnikov integral, we assumed that at $t=0$, the separatrix $\Gamma_0^{+}$ intersects section $J^c$ (see definition given in~\eqref{eq:DefGamma}), while
using the variational equations, at $t=0$, the separatrix intersects the section $J^{out}$.
Thus, let us take $s=s^{*}_1$ as the time it takes to go from $J^{out}$ to $J^c$ along the separatrix $\Gamma_0^{+}$.
Then, we replace $M(\theta)$ given in \eqref{eq:mel_int} by $\tilde M (\theta)= M (\theta + \omega s)$ in the map \eqref{eq:map_melnikov}.

Comparing the map \eqref{eq:map_melnikov} with $\tilde{M}(\theta)$ instead of $M(\theta)$ with expressions \eqref{eq:mapht} we have that the
following equalities must be satisfied:
\begin{equation}\label{eq:comparison}
\begin{array}{rcl}
\T^{*} &= & s^{*}_{0}+s^{*}_{1}, \\
\alpha u^{*}/v^{*} &= &1, \\
\eps \dosde u^{*} \rhoF(\theta)&=& \tilde{M}(\theta).\\
\end{array}
\end{equation}

\begin{figure}[h]
\begin{center}
\begin{tabular}{lll}
(a) &  (b) \\
\includegraphics[width=8cm]{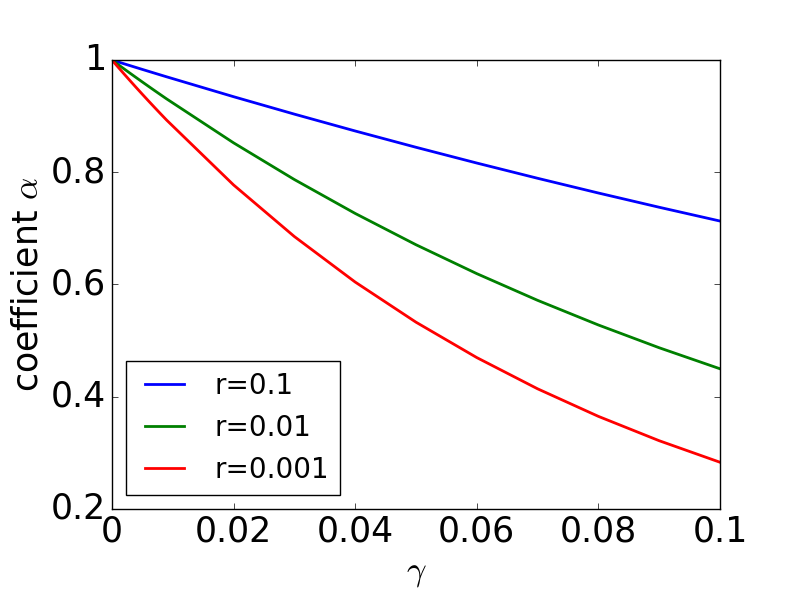} &
\includegraphics[width=8cm]{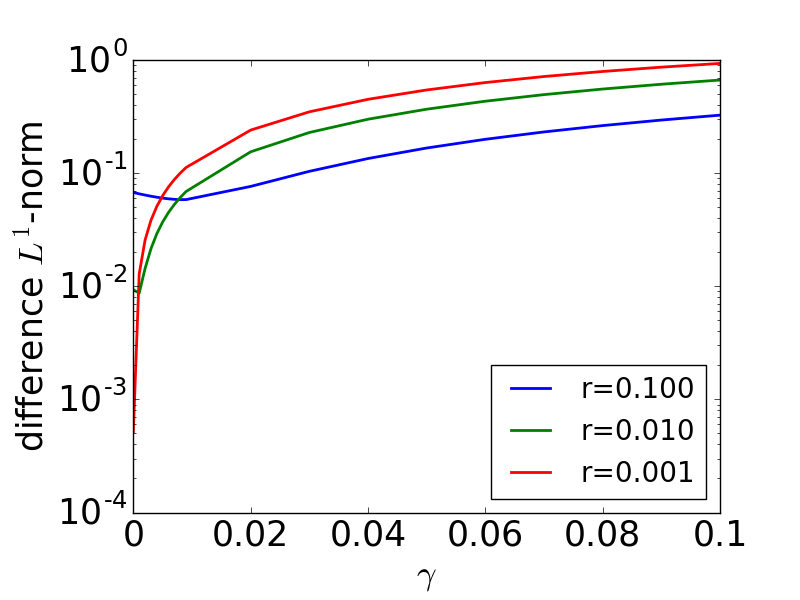} \\
\end{tabular}
\end{center}
\caption{
For different values of $\gamma$ and different values of $r$ of the Poincar\'{e} section we show
(a) the coefficient $\alpha$ in \eqref{eq:global_variacionals} computed numerically using variationals around
the perturbed separatrix and
(b) the difference in $L^1$-norm (in log-scale) between the function $\eps^{-1} \,\tilde M(\theta)$ and
the function $r \mu \rho(\theta)$ (see expression~\eqref{eq:comparison}).
}
\label{fig:compare_r}
\end{figure}

Next, we will show that these equalities are satisfied up to $\ocal(\gamma, u^{*},v^{*})$.
The first line can be checked straightforwardly and for the second line, since we have that
$\alpha=1+\ocal(\gamma)$ (see Figure~\ref{fig:compare_r}(a)), thus assuming that $u^{*}=v^{*}$, we have that there is agreement between both maps up to $\ocal(\gamma)$.

The third line must be checked numerically, and since we have chosen $\beta=5/4 \, \gamma$, $\tilde M$ is just
a trigonometric polynomial in $\theta$, so we need to compare the coefficients of the function
$\tilde{M}(\theta)$ with those of $\rhoF(\theta)$.
Indeed, $\rhoF(\theta)$ is a function of the form
$\rhoF(\theta)=\sum_i A_i \cos(\omega_i \theta) + B_i \sin (\omega_i \theta)$, where $A_i$ and $B_i$ satisfy $A_i=C_i \sin(\omega_i \phi_i)$ and
$B_i=C_i \cos (\omega_i \phi)$, with $\phi_i=1/\omega_i \arctan(A_i/B_i)$. Thus, we can write
\[ \rhoF (\theta) = \sum_{i=1}^n C_i \sin (\theta + \omega_i \phi_i)\]
where $C_i=A_i/\sin(\omega_i \phi_i)$, and clearly both functions $\rho(\theta)$ and $\tilde M(\theta)=M(\theta + \omega s)$
given in \eqref{eq:mel_int} have the same harmonics if $\phi_i=s$ for all $i$. This is true when $\gamma=0$
(results not shown) and for this case coefficients of $\mu u^{*} \rho(\theta)$ and $\eps^{-1} \tilde M (\theta)$ coincide up to an error which is $\ocal(r)$.
In Figure~\ref{fig:compare_r}(b) we show the comparison between both functions using the $L^1$-norm for different values of $\gamma$ and $r$. Clearly the error grows with $\gamma$ and $r$.

\section{Heteroclinic network model for binocular rivalry (HBR model)}\label{Sect:Heteroclinic}

In this section we want to apply the technique of the separatrix map explored in Section~\ref{Sect:Duffing} for the Duffing
equation to study models of bistable perception. We will consider a model of a specific phenomenon of bistable perception, namely,
binocular rivalry. In binocular rivalry two different images are presented to the two eyes simultaneously, and perception
alternates between these two images \cite{Huguet15}.
For this purpose we consider the model proposed by Ashwin and Lavric in \cite{Ashwin2010} that we will refer to as HBR model:

\begin{equation}\label{eq:ashwin_model}
\left\{
\begin{array}{rcl}
\dot{p} & = & h(p) + x^2 (1-p) + y^2 (-1-p), \\
\dot{x} & = & f(p,x,y) + I_x x + \eps \eta_x,\\
\dot{y} & = & g(p,x,y) + I_y y + \eps \eta_y,\\
\end{array}
\right.
\end{equation}
where $h(p) = -p (p-1)(p+1)$, $f(p,x,y)=((0.5-p)(p+1)-x^2-y^2)x$ and $g(p,x,y)=f(-p,y,x)$.
The variable $p$ represents the activity in the ``arbitration'' component, where $p=1$ ($p=-1$, respectively) represents
perception of the left (resp., right) eye stimulus, and $x$ and $y$ represent the activity pattern associated with stimulus
to the left and to the right eye, respectively. The quantities $I_x, I_y \geq 0$ represent external inputs to the system for $x$
and $y$, respectively.

The model has three equilibria for $\eps=0$, namely $(p,x,y)=(1,0,0)$ which corresponds
to the left dominant (LD) resting state, $(p,x,y)=(-1,0,0)$, which corresponds to the right dominant resting state (RD)
and $(p,x,y)=(0,0,0)$, which corresponds to a neutral state. The neutral state is asymptotically unstable with
eigenvalues $(1,0.5+I_x,0.5+I_y)$, while the LD and RD states are saddle points with eigenvalues $(-2, -1+I_x, I_y)$ and
$(-2,I_x,-1+I_y)$, respectively, for $I_x,I_y>0$. For the three equilibria, the associated eigenvectors are the canonical basis.
Moreover, $x=0$ and $y=0$ are invariant subspaces. The system has two heteroclinic orbits (by symmetry)
$\Gamma^{\pm}_{LD \rightarrow RD}$ from the LD to the RD saddle points lying on the plane $x=0$ and
two heteroclinic orbits (again by symmetry) $\Gamma^{\pm}_{RD \rightarrow LD}$
from the RD to the LD saddle points lying on the plane $y=0$ (see Figure~\ref{fig:SeparatricesAshwin}(a)).

As for the Duffing equation, we are going to consider a quasi-periodic perturbation representing external input to
the system for $x$ and $y$, i.e.
\begin{equation}\label{eq:per_per}
\eta_x (\theta)= \eta_y (\theta) = \sum_{i=1}^n a_i \cos(\theta_i), \quad \mbox{ with } \, \dot \theta=\omega,
\end{equation}
where $\theta=(\theta_1,\ldots,\theta_n) \in \mathbb{T}^n$, $\omega=(\omega_1,\ldots,\omega_n) \in \mathbb{R}^n$ and $a_i \in \mathbb{R}$ for $i=1,\ldots, n$.

In order to study the dynamics of system \eqref{eq:ashwin_model}-\eqref{eq:per_per} close to the heteroclinic cycles, we will compute the separatrix map.
As for the Duffing equation, in Section~\ref{sect:AshwinVar} we will provide details of the computation of the separatrix map
using variational equations along the heteroclinic cycles for $\eps=0$ and $I_x,\,I_y>0$.
In Section~\ref{sect:AshwinMel} we will provide details of the computation of the separatrix map via Melnikov integrals,
which for this model will be computed numerically.

Finally, we will compare the results with the case of a noisy perturbation,
\begin{equation}\label{eq:per_noiseni}
\eta_x=\eta_y=dW,
\end{equation}
where $dW$ is a zero mean Wiener process. We have also considered
$\eta_x = \eps_x\, dW_x$ and $\eta_y = \eps_y\, dW_y$, for two zero mean independent Wiener process
$dW_x$ and $dW_y$ as in \cite{Ashwin2010}, and numerical simulations (results not shown)
show qualitatively the same features.

\subsection{Derivation of the separatrix map via variational equations}\label{sect:AshwinVar}

\begin{figure}
\begin{center}
\begin{tabular}{lll}
(a) Global view &  (b) Partial view &  (c) Local view \\
\includegraphics[width=46truemm]{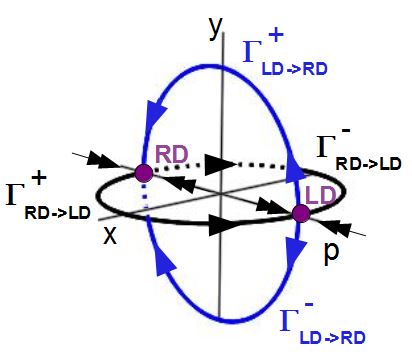}&
\includegraphics[width=46truemm]{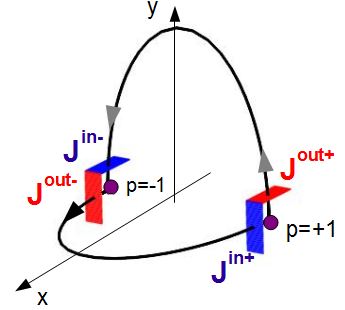} &
\includegraphics[width=46truemm]{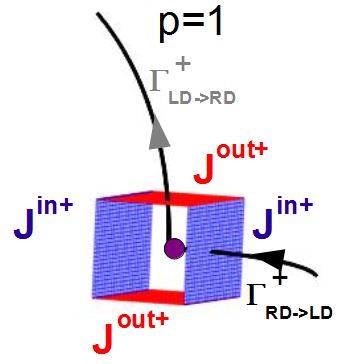}\\
\end{tabular}
\begin{tabular}{llll}
(d) & & & \\
\includegraphics[width=32truemm]{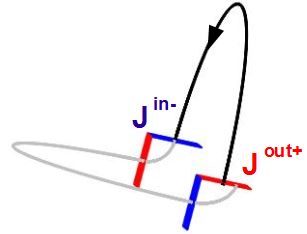} &
\includegraphics[width=30truemm]{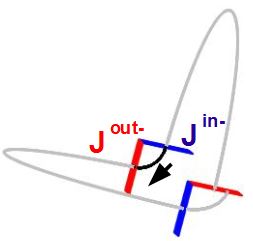} &
\includegraphics[width=30truemm]{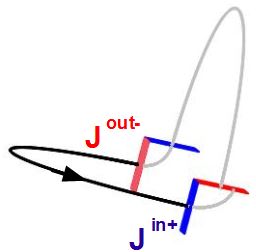} &
\includegraphics[width=30truemm]{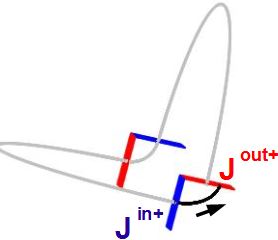} \\
$\Tg^{+}:\Sigma^{out+} \rightarrow \Sigma^{in-}$ &
$\Tl^{-}:\Sigma^{in-} \rightarrow \Sigma^{out-}$ &
$\Tg^{-}:\Sigma^{out-} \rightarrow \Sigma^{in+}$ &
$\Tl^{+}:\Sigma^{in+} \rightarrow \Sigma^{out+}$ \\
\end{tabular}
\end{center}
\caption{(a) Saddle points and heteroclinic orbits \eqref{eq:ashwin_model} for $I_x=I_y=\epsilon=0$.
(b, c) Views of the 2-D transversal sections used in the Poincar\'{e} maps defined in \eqref{eq:TransverseHN}.
Each Poincar\'{e} section, consists of two components, each one denoted with a sign subindex which is omitted here. (b) Positive components of the Poincar\'{e} sections (i.e., for $x>0$ and for $y>0$) both nearby $p=-1$ and $p=1$.
(c) All components of the Poincar\'{e} sections nearby $p=1$. Recall that $\Sigma^{in,out\, \pm}=J^{in,out\, \pm}\times \tor^n$, see \eqref{eq:ashwin_model}.
(d) The four Poincar\'{e} maps that define the total separatrix map $S=\Tl^{+} \circ \Tg^{-} \circ \Tl^{-} \circ \Tg^{+}$.}
\label{fig:SeparatricesAshwin}
\end{figure}

We consider system \eqref{eq:ashwin_model} with perturbation \eqref{eq:per_per} and compute the return map using variational equations
as in Section \ref{sect:VariationalDuffing}. Here we have one more state variable, say $p$, but the procedure is analogous.
First, we define 2-dimensional sections close to the LD saddle point with $p=1$ and to the RD saddle point with $p=-1$, both transversal to the flow of
system \eqref{eq:ashwin_model} for $\eps=0$ (see also Figure~\ref{fig:SeparatricesAshwin}(b) and (c)):
\begin{equation}\label{eq:TransverseHN}
\begin{array}{rcl}
J^{out+}_{\sigma} & = & \{(p=1+q, x, \sigma y^{*})\}, \\
J^{in+}_{\sigma} & = & \{(p=1+q, \sigma x^{*}, y)\}, \\
J^{out-}_{\sigma} & = & \{(p=-1+q, \sigma x^{*}, y)\}, \\
J^{in-}_{\sigma} & = & \{(p=-1+q, x, \sigma y^{*})\}, \\
\end{array}
\end{equation}
where $\sigma \in \{ +, -\}$ and $x^{*},y^{*}>0$ are fixed. When considering the angular variable $\theta \in \mathbb{T}^n$, we denote by
\begin{equation}\label{eq:sigma_sec}
\Sigma_{\sigma}^{in\pm,out\pm}=J_{\sigma}^{in\pm,out\pm}\times \mathbb{T}^n
\end{equation}
the corresponding sections in the extended phase space.

Notice that for the Duffing equation we included a variable $\sigma$ that allows us to distinguish whether the trajectory hits the
section $\Sigma^{out \pm,in \pm}_{+}$ or $\Sigma^{out \pm,in \pm}_{-}$. However, since the system has the symmetry
$(p,x,y) \rightarrow (-p,y,x)$, the information provided by $\sigma$ does not affect the results, so we
are not going to include $\sigma$ to simplify the reading. For the same reason, we will denote
\begin{equation}\label{def:SigmaAshwin}
\Sigma^{in \pm,out\pm}=\Sigma^{in \pm,out\pm}_+ \cup \Sigma^{in \pm,out\pm}_{-},
\end{equation}
and only consider two heteroclinic connections, namely $\Gamma_{RD \rightarrow LD}^{+}$ and $\Gamma_{LD \rightarrow RD}^{+}$,
for which we will omit the superindex from now on.

First, we are going to obtain an approximation of the local maps $\Tl^\pm$ from section $\Sigma^{in \pm}$ to $\Sigma^{out \pm}$
using the linearized dynamics around the saddle points with $p=1$ for $\Tl^+$ and $p=-1$ for $\Tl^{-}$ (see Figure~\ref{fig:SeparatricesAshwin}(d)). Thus, we have
\begin{equation}\label{eq:ashwin_local}
\begin{array}{rcl}
\Tl^{+}:& \Sigma^ {in+} & \longrightarrow \Sigma^{out+}\\
& (q,y,\theta) & \mapsto (\bar{q},\bar{x},\bar{\theta}),
\end{array}
\end{equation}
where
\begin{equation}\label{eq:local_ashwin1}
\begin{array}{rl}
\bar{q}&=q \exp\left(-2\, \T_y\right),\\
\bar{x}&=x^{*}\,\exp\left((-1+I_x)\, \T_y\right),\\
\bar{\theta}&=\theta+\omega\, \T_y,\\
\end{array}
\end{equation}
with
\begin{equation}\label{eq:time_ashwin_local}
\T_y=\displaystyle\frac{1}{I_y} \ln \left ( \frac{y^{*}}{y}\right).
\end{equation}

Similarly, the map $\Tl^{-}: \Sigma^{in-} \longrightarrow \Sigma^{out-}$ is described by
\begin{equation}\label{eq:local_ashwin2}
(\bar{q},\bar{y},\bar{\theta})=\Tl^{-}(q,x,\theta)=(q \exp\left(-2\,\T_x\right),y^{*}\, \exp\left((-1+I_y)\,\T_x\right),\theta+\omega\, \T_x),
\end{equation}
with $\T_x=\displaystyle\frac{1}{I_x} \ln \left ( \frac{x^{*}}{x}\right)$.

Next, we are going to obtain an approximation of the global maps $\Tg^{\pm}$ from sections $\Sigma^{out \pm}$ to $\Sigma^{in \mp}$
(see Figure~\ref{fig:SeparatricesAshwin}(d)). Consider first the global map
\begin{equation}\label{eq:ashwin_global}
\begin{array}{rcl}
\Tg^{+}:& \Sigma^{out+} & \longrightarrow \Sigma^{in-}\\
&(q,x,\theta) & \mapsto (\bar{q},\bar{x},\bar{\theta}),
\end{array}
\end{equation}
defined as
\[\Tg^+ (q,x,\theta) = \varphi (\tau^{*}(1+q,x,y^{*},\theta);1+q,x,y^{*},\theta),\]
where $\varphi$ is the flow of system \eqref{eq:ashwin_model}-\eqref{eq:per_per} and
$\tau^{*}(1+q,x,y^{*},\theta)>0$ is the minimal global time such that $\varphi(\tau^{*}(1+q,x,y^{*},\theta);1+q,x,y^{*},\theta) \in \Sigma^{in-}$.

We will obtain an approximation of the global map by computing the linear dynamics around
the separatrix $\Gamma_{LD \rightarrow RD}$  for $I_x,I_y\in (0,1)$ and $\eps=0$. Thus,
as in Section~\ref{sect:VariationalDuffing}, we consider the parameter $\eps$ as a variable ($\dot{\eps}=0$) and denote by $\hat \varphi (t;p,x,y,\theta,\eps)$ the
flow of the extended system \eqref{eq:ashwin_model}-\eqref{eq:per_per} adding
$\dot \eps = 0$, and $\hat T_G^{+}$ the extended global map.
Consider a point of the form
$(1+q,x,y^{*},\theta,\eps)=(1 + q^s + \Delta q, x^s + \Delta x, y^{*}, \theta^s + \Delta \theta, \eps)$, where
$\{(1+q^s,x^s,y^{*})\}= \Gamma_{LD \rightarrow RD} \cap J^{out+} $ and $\theta^s \in \mathbb{T}^n$.
Notice that $x^s=0$.  Its image is given by
\[
\begin{array}{rcl}
\hat T_G^{+}(q,x,y^{*},\theta,\eps) & = & \hat T_G^{+}(q^s,x^s,y^{*},\theta^s,0)+D \hat T_G^{+}(q^s,x^s,y^{*},\theta^s,0) \cdot \Delta + \ocal(\Delta^2) \\
& = &  (\bar{q}^s,\bar{x}^s,y^{*},\theta^s + \omega \T^{*},0) \\
& & + (\alpha_q \Delta x + \beta_q \Delta q + \eps \rhoF_q(\theta^s),
\alpha_x \Delta x + \beta_x \Delta q + \eps \rhoF_x(\theta^s),\\
& & \qquad
0,\Delta \theta + \alpha_{\theta} \Delta x + \beta_{\theta} \Delta q +\eps \rhoF_{\theta}(\theta^s),
\eps) \\
& & + \ocal(\Delta^2)\\
\end{array}
\]
where $\Delta=(\Delta q, \Delta x, 0, \Delta \theta, \eps)$, $\{(-1+\bar{q}^s,\bar{x}^s,y^{*})\}= \Gamma_{LD \rightarrow RD} \cap J^{in-}$, $\T^{*}=\tau^{*}(1+\bar{q}^s,x^s,y^{*},\theta^s)$, which is independent of $\theta^s$ when $\eps=0$.
Moreover, analogously to
\eqref{eq:coef_uv} but taking into account that now there is an extra variable $p$, we have the following expression for the coefficients:
\begin{equation}\label{eq:coef_ashwin_fun}
\begin{array}{lll}
\alpha_x=\hp^x_x - \dfrac{F^x}{F^y} \hp^y_x, & \beta_x=\hp^x_p - \dfrac{F^x}{F^y}\hp^y_p, & \rhoF_x(\theta)=\hp^x_{\eps} - \dfrac{F^x}{F^y} \hp^y_{\eps}, \vspace{0.2cm} \\
\alpha_q=\hp^p_x - \dfrac{F^p}{F^y} \hp^y_x, & \beta_q=\hp^p_p - \dfrac{F^p}{F^y} \hp^y_p, &  \rhoF_q(\theta)=\hp^p_{\eps} - \dfrac{F^p}{F^y} \hp^y_{\eps}, \vspace{0.2cm} \\
\alpha_{\theta}=-\omega \dfrac{\hp^y_x}{F^y}, & \beta_{\theta}=-\omega \dfrac{\hp^y_p}{F^y}, &  \rhoF_{\theta}(\theta)=- \omega \dfrac{\hp^y_{\eps}}{F^y}, \vspace{0.2cm} \\
\end{array}
\end{equation}
where $F^w$ and $\hat \varphi^w$ denote the $w$-coordinate of the vector field \eqref{eq:ashwin_model} and the
extended flow $\hat \varphi$, respectively, for $w=p,x,y$.

Therefore, disregarding the terms of $\ocal(\Delta^2)$, the global map \eqref{eq:ashwin_global} has the form below, which is
analogous to expression~\eqref{eq:exp_global} taking into account that there is an extra variable $p$ and $x^s=0$:
\begin{equation}\label{eq:global_variacionals_ashwin}
\begin{array}{rcl}
\bar{q} &= &\bar{q}^s + \alpha_q\, (x - x^s)  + \beta_q\, (q - q^s)+ \eps\, \rhoF_q(\theta), \\
\bar{x} &= & \alpha_x\, (x - x^s) + \beta_x\, (q - q^s) + \eps\, \rhoF_x(\theta), \\
\bar{\theta}&=& \theta+ \omega\, \T^{*}+ \alpha_{\theta}\, (x - x^s) + \beta_{\theta}\, (q - q^s) +\eps\, \rhoF_{\theta}(\theta).\\
\end{array}
\end{equation}

Since $x-x^s$, $q-q^s$ and $\eps$ are assumed to be small, the contribution of the terms $\alpha_{\theta}\, (x - x^s)$,
$\beta_{\theta}\, (q - q^s)$ and $\eps\, \rhoF_{\theta}$ is neglegible compared with the
finite term $\omega \T^{*}$. Recall $x^s=0$. Moreover, one can see that the terms $\bar{q}^s$ and $q^s$ are $\ocal_2(y^{*},x^{*})$.
Therefore, considering only the dominant terms, we can write the global map $\Tg^+$ as
\begin{equation}\label{eq:global_variacionals_ashwin_simplified}
\begin{array}{rcl}
\bar{q} &= & \alpha_q\, x + \beta_q\, q + \eps\, \rhoF_q(\theta), \\
\bar{x} &= & \alpha_x\, x + \beta_x \, q + \eps\, \rhoF_x(\theta), \\
\bar{\theta}&=& \theta + \omega\, \T^{*}. \\
\end{array}
\end{equation}

The global map
$$
\begin{array}{rl}
\Tg^{-}: \Sigma^{out-} & \longrightarrow \Sigma^{in+}\\
(q,y,\theta) & \mapsto (\bar{q},\bar{y},\bar{\theta}),
\end{array}
$$
is the same as \eqref{eq:global_variacionals_ashwin_simplified} because of the symmetry $(p,x,y) \mapsto (-p,y,x)$,
just replacing $x$ by $y$ and $\bar{x}$ by $\bar{y}$.
that is,
\begin{equation}\label{eq:global_variacionals_ashwin2}
\begin{array}{rcl}
\bar{q} &= & \alpha_q\, y + \beta_q\, q + \eps\, \rhoF_q(\theta), \\
\bar{y} &= & \alpha_x\, y + \beta_x\, q + \eps\, \rhoF_x(\theta), \\
\bar{\theta}&=& \theta + \omega\, \T^{*}. \\
\end{array}
\end{equation}

We finally define the separatrix map $S$ (as in Section~\ref{sec:ss_ssve}) by combining the local and global maps described above in the following way
(see Figure~\ref{fig:SeparatricesAshwin}(d)):
\begin{equation} \label{eq:ashwin_map}
\begin{array}{rccc}
S:=\Tl^{+} \circ \Tg^{-} \circ \Tl^{-} \circ \Tg^{+}:& \Sigma^{out+} &\rightarrow & \Sigma^{out+} \\
& (q,x,\theta) & \rightarrow & (\hat{q},\hat{x},\hat{\theta}). \\
\end{array}
\end{equation}
For the sake of clarity we will not provide the explicit global expression of the map $S$, and work only with the expressions
\eqref{eq:local_ashwin1}-\eqref{eq:local_ashwin2}-\eqref{eq:global_variacionals_ashwin_simplified}-\eqref{eq:global_variacionals_ashwin2}.

We conclude this section with an important remark. For this particular model
$\alpha_q \approx 0$, $\beta_q \approx 0$, $\alpha_{\theta} \approx 0$ and $\beta_x \approx 0$ in expressions
\eqref{eq:global_variacionals_ashwin_simplified} and \eqref{eq:global_variacionals_ashwin2} (see Section~\ref{Sect:AshwinNumeric}).
Therefore, by looking at these expressions and
the ones for the local maps in \eqref{eq:local_ashwin1}-\eqref{eq:local_ashwin2}, it is clear that the dynamics of the variable $q$ decouples from the
rest of the variables. Therefore, \eqref{eq:ashwin_map} will be essentially a 2-dimensional map, since we can omit the dynamics for $q$
(see also \cite{Ashwin2010}).

\subsection{Derivation of the separatrix map via Melnikov integrals (semi-analytical derivation)}\label{sect:AshwinMel}

We can also obtain the separatrix map via Melnikov integrals, as done in the study of the Duffing equation, see Section \ref{sect:MelnikovDuffing}.
In the case of system \eqref{eq:ashwin_model} with $I_x=I_y=\eps=0$, we can take advantage of the knowledge of the heteroclinic connections
given by $\NansaGx0=\{\varG=0\} \cap \{x=0\}$ and $\NansaFy0=\{\varF=0\} \cap \{y=0\}$, where $\varG:=p^2+2\,y^2-1$ and $\varF:=p^2+2\,x^2-1$.
Let us denote the respective time-parameterizations by
$\{(p_0(s),0,y_0(s)), s\in\rr\}$, where $\lim_{s\to \mp \infty}p_0(s)=\pm 1$ and $\lim_{s\to \mp \infty}y_0(s)=0$, and
$\{(p_0(-s),x_0(s),0), s\in\rr\}$, where $\lim_{s\to \mp \infty}x_0(s)=0$. More specifically, $p_0(s)=\xi\left(1+\xi^2\right)^{-1/2}$
and $x_0(s)=y_0(s)=\left(2\,(1+\xi^2)\right)^{-1/2}$, where $\xi=\xi(s)$ is implicitly given by
$s=\xi\,\sqrt{\xi^2+1}+\mbox{arcsinh}(\xi)-\xi^2$.
The variables $\varG$ and $\varF$ play an analogous role to the Hamiltonian variable
$h$ used in Section \ref{sect:MelnikovDuffing}.

We use the Poincar\'{e} sections already defined in \eqref{eq:sigma_sec}. For the sections $\Sigma^{out+}$ and $\Sigma^{in-}$ we will consider the set of variables $(\varG,x,\theta)$, while for the sections $\Sigma^{in+}$ and $\Sigma^{out-}$ we will consider
the set of variables $(\varF,y,\theta)$. As in Section \ref{sect:AshwinVar}, the separatrix map is composed by the concatenation
of the four Poincar\'{e} maps illustrated in Figure~\ref{fig:SeparatricesAshwin}(d), that is: $S=\Tl^{+} \circ \Tg^{-} \circ \Tl^{-} \circ \Tg^{+}$.

The local maps are defined as in \eqref{eq:ashwin_local} and we only need to rewrite them in the new variables. For instance,
$$
\begin{array}{rccl}
\Tl^{+}:& \Sigma^{in+} & \rightarrow & \Sigma^{out+}\\
& (\varF,y,\theta) & \mapsto & (\varGbar,\bar{x},\bar{\theta}),
\end{array}
$$
where
$$
\begin{array}{rl}
\varGbar&=\left(\varF-2\,{x^{*}}^2\right)\,\psi^2+2\,{y^{*}}^2+2\,\psi\,(1-\psi)\,\left(-1+\sqrt{1+\varF-2\,{x^{*}}^2}\right),\\
\bar{x}&=x^{*}\exp\left((-1+I_x)\, \T_y\right),\\
\bar{\theta}&=\theta+\omega\,\T_y,
\end{array}
$$
where the time $\T_y$ is defined in \eqref{eq:time_ashwin_local} and
$\psi:=\exp\left(-2\,\T_y\right)$. The map $\Tl^{-}: \Sigma^{in-} \longrightarrow \Sigma^{out-}$
is obtained in a similar way and will be described as $(\varFbar,\bar{y},\bar{\theta})=\Tl^{-}(\varG,x,\theta)$.

To compute the global maps
$$
\begin{array}{rcclcrccl}
\Tg^{+}:& \Sigma^{out+} & \rightarrow & \Sigma^{in-}&\qquad&\Tg^{-}:& \Sigma^{out-} & \rightarrow & \Sigma^{in+}\\
& (\varG,x,\theta) & \mapsto & (\varGbar,\bar{x},\bar{\theta}),&\qquad& & (\varF,y,\theta)& \mapsto & (\varFbar,\bar{y},\bar{\theta}),
\end{array}
$$
we use a different strategy. For instance, to get $\Tg^{+}$ up to the first order in terms of the perturbation parameters $I_x$, $I_y$ and $\eps$ (similar expressions can be obtained for $\Tg^{-}$), we remark that $\varG=x=0$ (the heteroclinic connection $\NansaGx0$) is invariant
for the flow of system \eqref{eq:ashwin_model}-\eqref{eq:per_per} when $I_x=I_y=\eps=0$. Therefore, $\varG$ and $x$ satisfy differential equations that can be written as
\begin{equation}\label{eq:VarsGx}
\left(
\begin{array}{c}
\dotvarG\\
\dot{x}
\end{array}
\right)=
\left(
\begin{array}{cc}
a & b\\
0 & c
\end{array}
\right)\,
\left(
\begin{array}{c}
\varG\\
x
\end{array}
\right)
+
\left(
\begin{array}{c}
m\\
n
\end{array}
\right)
,
\end{equation}
where $m$, $n$ vanish when $I_x=I_y=\eps=0$. More specifically, $a,b,c,m,n$ are given by:
$$
\begin{array}{rl}
a&=a(p,x,y):=-2(p^2+x^2+y^2),\\
b&=b(p,x):=-2x^2(1-p),\\
c&=c(p,x,y):=(0.5-p)(p+1)-x^2-y^2,\\
m&=m_I(y)\,I_y+m_{\eps}(\theta)\,\eps:=(4\,y^2)\,I_y+(4 \,\eta_y(\theta) y)\,\eps,\\
n&=n_I(x)\,I_x+n_{\eps}(\theta)\,\eps:=x\,I_x+\eta_x(\theta)\,\eps.
\end{array}
$$
Note that if we approximate system \eqref{eq:VarsGx} by evaluating the functions $a$, $b$, $c$, $m$ and $n$ on $\NansaGx0$,
that is, with $p=p_0(s)$, $x=0$, $y=y_0(s)$, then \eqref{eq:VarsGx} becomes an uncoupled system of linear ordinary differential equations in $\varG$ and $x$:
\begin{equation}\label{eq:VarsGx0}
\left(
\begin{array}{c}
\dotvarG\\
\dot{x}
\end{array}
\right)=
\left(
\begin{array}{cc}
a_0(s) & 0\\
0 & c_0(s)
\end{array}
\right)
\left(
\begin{array}{c}
\varG\\
x
\end{array}
\right)
+
\left(
\begin{array}{c}
m_0(s)\\
n_0(s)
\end{array}
\right),
\end{equation}
where
$a_0(s)=-2(p_0(s)^2+y_0(s)^2)$, $m_0(s)=(4\,y_0(s)^2)\,I_y+(4 \,\eta_y(\theta+\omega\,s) y)\,\eps$,
$c_0(s)=(0.5-p_0(s))(p_0(s)+1)-y_0(s)^2$ and $n_0(s)=\eta_x(\theta+\omega\,s)\,\eps$.

We introduce $\tau$ and $\bar{\tau}$ such that $(p_0(\tau),0,y_0(\tau))\in J^{out+}$ and $(p_0(\bar{\tau}),0,y_0(\bar{\tau}))\in J^{in-}$ and set
$\varG=\varG(\tau)$, $\varGbar=\varG(\bar{\tau})$, $x=x(\tau)$ and $\bar{x}=x(\bar{\tau})$. Then, we can solve the two uncoupled differential equations
\eqref{eq:VarsGx0} and obtain an expression for the Poincar\'{e} map $\Tg^{+}$ up to order one in terms of the perturbation parameters $I_x$, $I_y$ and $\eps$, which can be seen as analogous to the Melnikov integrals
used in \eqref{eq:DuffingGlobalMap} for the Duffing equation. More precisely, the map has the form
\begin{equation}\label{eq:AshwinMapMelnikov}
\begin{array}{rl}
\varGbar & = B_{h}\, \varG + \eps \, P_{h}(\theta) + \Xi_{I} \, I_y , \\
\bar{x} & = A_x\, x + \eps\, P_{x}(\theta), \\
\bar{\theta} &= \theta + \omega\, (\bar{\tau} - \tau), \\
\end{array}
\end{equation}
where
\begin{equation}\label{eq:param_map}
\begin{array}{l}
 B_{h}= \exp \left (\int_{\tau}^{\bar{\tau}} a_0(s) ds \right),\\
 P_{h}(\theta)=\int_{\tau}^{\bar{\tau}} 4 y_0(t) \eta_y(\theta+\omega\,t) \exp\left(-\int_{\tau}^t a_0(s) ds \right) dt,\\
 \Xi_{I}=\int_{\tau}^{\bar{\tau}} 4 y_0^2(t) \exp \left (-\int_{\tau}^t a_0(s) ds \right ) dt,\\
 A_x= \exp \left ( \int_{\tau}^{\bar{\tau}} c_0(s) ds \right),\\
 P_{x}(\theta)=\int_{\tau}^{\bar{\tau}} \eta_x(\theta+\omega\,t) \exp \left ( - \int_{\tau}^t c_0(s) ds \right ) dt.
\end{array}
\end{equation}
The expression \eqref{eq:AshwinMapMelnikov}-\eqref{eq:param_map} provides a theoretical framework to compute an approximation
of the global map $\Tg^{+}$ for $I_x$, $I_y$ and $\epsilon$ small enough, but the cumbersome parameterization of $\NansaGx0$
makes a numerical resolution more advisable (see Section~\ref{Sect:AshwinNumeric}).

\subsection{Numerical computation}\label{Sect:AshwinNumeric}

In this Section, we show the numerical computations of the separatrix map for the HBR model with perturbation \eqref{eq:per_per} and the same frequencies as in Section~\ref{sec:numerical}, that
is $\omega=(\omega_1,\omega_2,\omega_3)=(1,(\sqrt{5}-1)/2,\sqrt{769}-27)$. Moreover, we
choose sections $\Sigma^{out \pm, in \pm}$ with $x^{*}=y^{*}=r$ (see definition~\eqref{eq:TransverseHN} for the sections)
and $r=0.1$ (see Appendix~\ref{ap:cr} for a justification for this choice). Moreover, we consider inputs
$I=I_x=I_y=0.1$ as in \cite{Ashwin2010}.

For these parameter values we have computed the global map \eqref{eq:global_variacionals_ashwin} obtained using variational equations.
We have also computed the global map \eqref{eq:AshwinMapMelnikov} using the alternative method based on Melnikov integrals, and the comparison between both methods is
discussed in Section~\ref{sec:com_het}. For the rest of the section, we focus on the map computed using variational equations, since it
is applicable to a wider parametric range of $I$ values.

For the specific choice of the parameters described above, we compute numerically the coefficients of the map \eqref{eq:global_variacionals_ashwin} using the expressions in \eqref{eq:coef_ashwin_fun}, and we obtain:
\begin{equation}\label{eq:map_ashwin_calculat}
\begin{array}{rcl}
\bar{q}&=& 0.005494470100 + \eps\, \rhoF_q(\theta), \\
\bar{x} & = & 0.0000123595\, x + \eps\, \rhoF_x(\theta), \\
\bar{\theta}&=& \theta + 19.2385452050 \, \omega + 7.1811476867\, (q+0.0091291201)+\eps\, \rhoF_{\theta}(\theta), \\
\end{array}
\end{equation}
where
\[
\begin{array}{rl}
\rhoF_q(\theta)=& a_1 ( -0.0830186087 \cos(\theta_1)-0.0388801779\sin(\theta_1) )\\
&+a_2(-0.0355217259\cos(\theta_2) -0.1058303244\sin(\theta_2)) \\
&+a_3(-0.0708199918 \cos(\theta_3)+0.0786695423 \sin(\theta_3)),
\end{array}
\]
\[
\begin{array}{rl}
\rhoF_x(\theta)=&a_1(-0.4340559240 \cos(\theta_1)+0.7770758314 \sin(\theta_1))\\
&+a_2(-2.9264485016 \cos(\theta_2)+ 1.8586408166 \sin(\theta_2)) \\
&+a_3(1.9947756545 \cos(\theta_3)+1.5403924072 \sin(\theta_3)),
\end{array}
\]
and
\[
\begin{array}{rl}
\rhoF_{\theta}(\theta)=& a_1(-3.9499622400 \cos(\theta_1)+97.4619562829 \sin(\theta_1)) \\
&+a_2(-28.9752119958 \cos(\theta_2)+158.9382104783 \sin(\theta_2)) \\
&+a_3(-13.2998112505 \cos(\theta_3)+123.3597165763 \sin(\theta_3)).
\end{array}
\]
Notice that in expression \eqref{eq:map_ashwin_calculat} we have $\alpha_q=\beta_q=\beta_x=\alpha_{\theta}=0$.
Its reduced version \eqref{eq:global_variacionals_ashwin_simplified} reads out as
\begin{equation}\label{eq:map_ashwin_reduced_comp}
\begin{array}{rcl}
\bar{q} &= & \eps\, \rhoF_q(\theta), \\
\bar{x} &= &0.0000123595\, x + \eps\, \rhoF_x(\theta), \\
\bar{\theta}&=& \theta + 19.2385452050 \, \omega.
\end{array}
\end{equation}

The local maps are obtained explicitly according to the formulas in \eqref{eq:local_ashwin1}-\eqref{eq:local_ashwin2}.
Thus, we consider the separatrix map
$S$ defined in \eqref{eq:ashwin_map} using the specific parameters computed in \eqref{eq:map_ashwin_reduced_comp}
for the reduced version of the global map. Next, we are going to explore the dynamics of this map (omitting the dynamics for $q$ which decouples from the other two variables)
for $\eps=0.001$ as in \cite{Ashwin2010} and
three different perturbations: a periodic perturbation ($a_1=1, a_2=0, a_3=0$), a quasi-periodic perturbation with 2 frequencies ($a_1=1,a_2=1,a_3=0$) and
a quasi-periodic perturbation with 3 frequencies ($a_1=1,a_2=1,a_3=1$).

\begin{figure}
\begin{center}
\begin{tabular}{lll}
(a) 1 frequency & (b) 2 frequencies &  (c) 3 frequencies \\
\includegraphics[width=5.5cm]{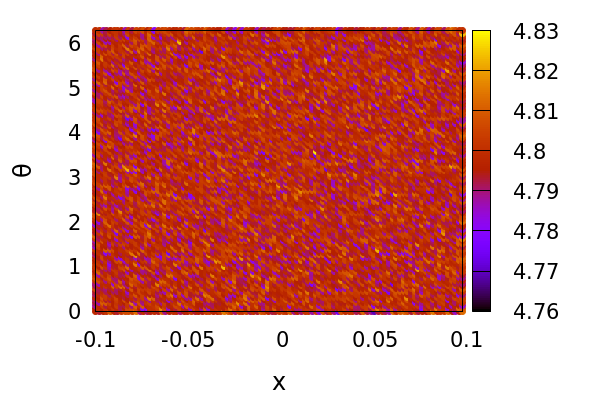} &
\includegraphics[width=5.5cm]{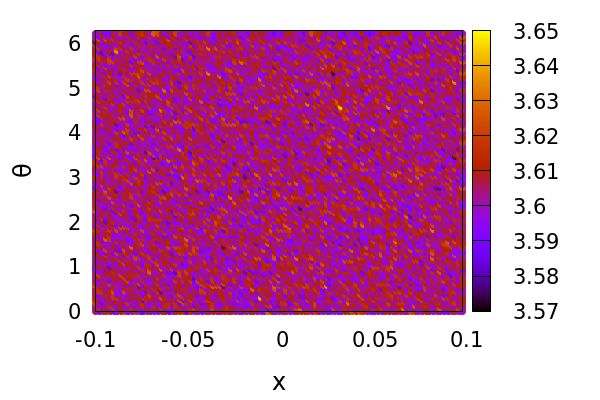} &
\includegraphics[width=5.5cm]{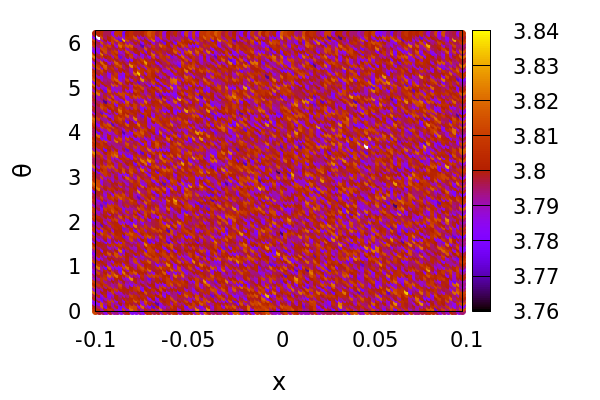} \\
\end{tabular}
\caption{Maximal Lyapunov exponent (computed using MEGNO) for the orbits of the separatrix map \eqref{eq:ashwin_map} (using the reduced version \eqref{eq:map_ashwin_reduced_comp} of the global map) for the
HBR model with initial conditions on the phase space $(x,\theta)$ ($\theta_i=\theta$, for $i=1,2,3$) and parameters $I_x=I_y=0.1$, $\eps=0.001$ and $r=0.1$.
The number of frequencies in the perturbation is indicated in each panel.}
\label{fig:megno_ashwin}
\end{center}
\end{figure}

To explore the dynamics for each map we consider a grid of initial conditions on the plane $(x,\theta)$, where all the components of
the vector $\theta$ take the same value $\theta_1=\theta_2=\theta_3$, and we compute for each
corresponding orbit the maximal Lyapunov exponent using the MEGNO algorithm.
Results are shown in Figure~\ref{fig:megno_ashwin}. We observe that,
independently of the initial conditions, all the orbits show chaotic behaviour (positive
Lyapunov exponent) for the three maps (1,2 or 3 frequencies).

\begin{figure}[h]
\begin{center}
\begin{tabular}{ll}
(a) 1 frequency & (b) 2 frequencies \\
\includegraphics[width=6cm]{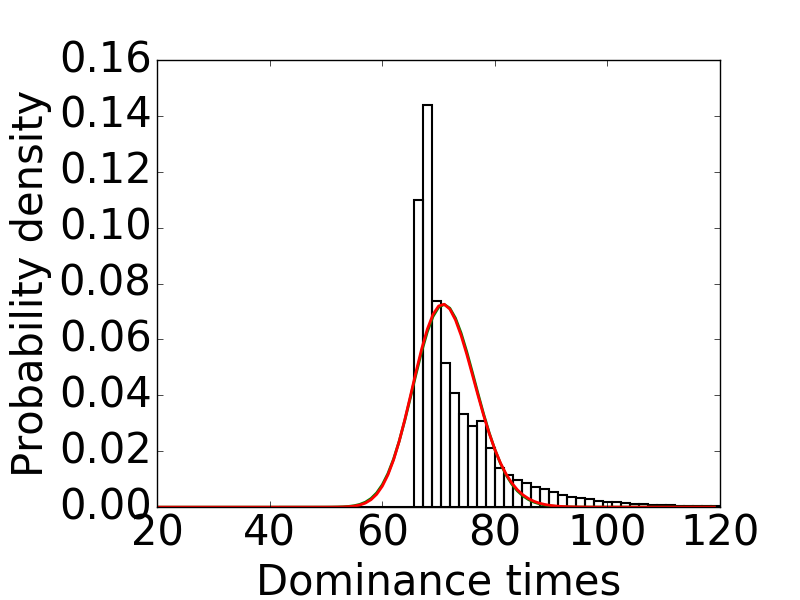}&
\includegraphics[width=6cm]{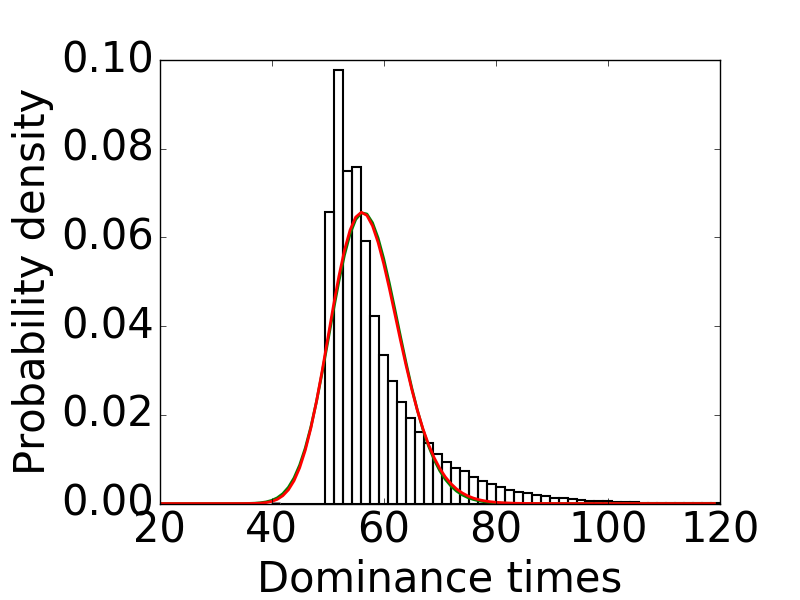}\\
(c) 3 frequencies & (d) noise \\
\includegraphics[width=6cm]{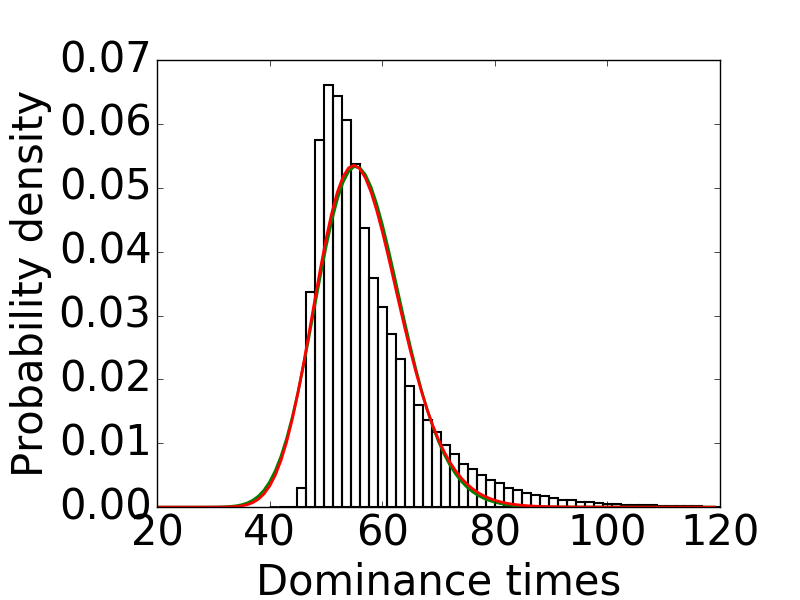} &
\includegraphics[width=6cm]{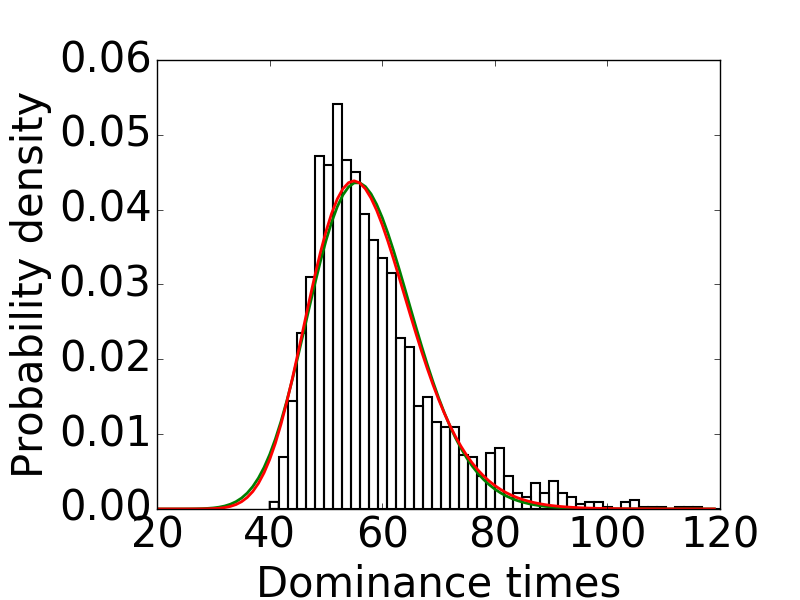}\\
\end{tabular}
\caption{
(a,b,c) Dominance times histograms computed using the separatrix map \eqref{eq:ashwin_map} with a quasi-periodic perturbation of 1,2 and 3
frequencies \eqref{eq:per_per} and $\eps=10^{-3}$ and $I_x=I_y=0.1$.
We show time to return to the section $p=0$ by concatenating one local map and one global map.
We used 200,000 iterates. Initial conditon is $x=-0.1$, $\theta_i=0$.
(d) Dominance times for system \eqref{eq:ashwin_model} with noise \eqref{eq:per_noiseni} with  $\eps=10^{-3}$ and $I_x=I_y=0.1$.
We show time to return to the Poincar\'{e} section $p=0$ by integrating the
full system \eqref{eq:ashwin_model} (notice that this is equivalent to concatenate one local map and one global map when we consider symmetry in the system).
We used 2,000 iterates. For all histograms, we show fittings to Gamma (green) and log-normal (red) distributions (fittings not distinguishable). Histograms have been normalized to have area 1.
}
\label{fig:DominanceAshwin}
\end{center}
\end{figure}

\begin{figure}[h]
\begin{center}
\includegraphics[width=7cm]{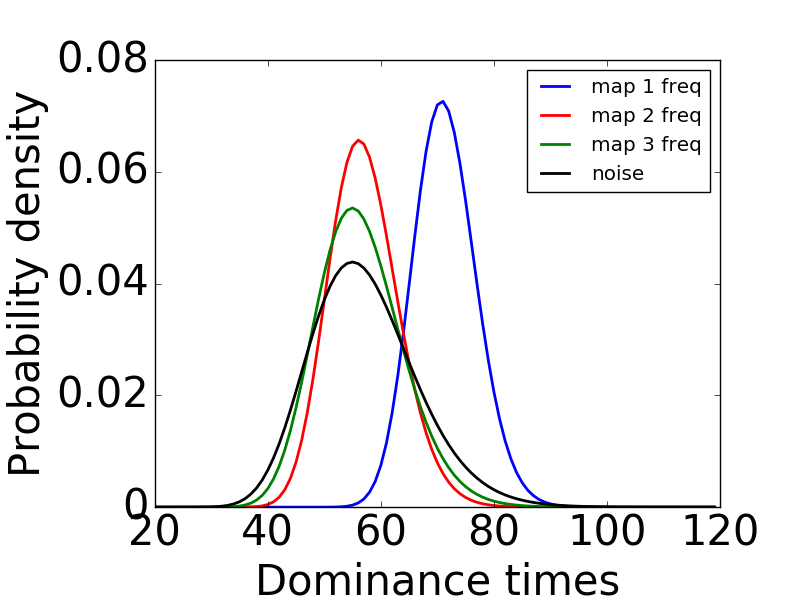}
\end{center}
\caption{Comparison between fittings to log-normal distribution of the dominance times histograms in Figure~\ref{fig:DominanceAshwin}  (similar results for fittings to Gamma distribution).}
\label{fig:resum_hist}
\end{figure}

Moreover, we explored the distribution of dominance times $\T_{dom}$
defined as the time difference between impacts on the sections $\Sigma^{out+}$ and $\Sigma^{out-}$, i.e
$\T_{dom}=(\tilde \theta_i-\theta_i)/\omega_i$, where $\tilde \theta_i$ is the projection onto the $\theta_i$ component of $T_L \circ T_G$.
Notice that $\T_{dom}$ is independent of the coordinate $i$.
The dominance times provide an approximation of the time spent near the vicinity of the saddle points
$p=1$ (corresponding to the LD state) and $p=-1$ (corresponding to the RD state).
Notice that this is equivalent to consider
the time difference between impacts on the section $\Sigma^c$ corresponding to $p=0$.
The histograms of $\T_{dom}$ for 1, 2 and 3 frequencies are shown in Figure~\ref{fig:DominanceAshwin}(a), (b) and (c), respectively.
We compare these histograms with those obtained with a noisy
perturbation~\eqref{eq:ashwin_model}-\eqref{eq:per_per} (Figure~\ref{fig:DominanceAshwin}(d)).
We fit the histograms to Gamma and log-normal distributions (see Appendix~\ref{ap:fitting}). Notice that differences are not noticeable.
In Figure~\ref{fig:resum_hist} we show the fittings altogether. We observe that as the number of frequencies in the perturbation increases the histograms
of the dominance durations shift leftwards and they become more similar to the ones
obtained with noise.

\subsubsection{Comparison between separatrix maps for the HBR model}\label{sec:com_het}

We have performed numerical simulations and computed the coefficients in \eqref{eq:param_map}
for $x^{*}=y^{*}=r=0.1$ and the same three frequencies as before, that is $\omega=(\omega_1,\omega_2,\omega_3)=(1,(\sqrt{5}-1)/2,\sqrt{769}-27)$, which lead to:
\begin{equation}\label{eq:AshwinMapMelnikovExample}
\begin{array}{rl}
 B_{h} = & 0,\\
 P_{h}(\theta)= &
  a_1 (-0.2774816480\,\cos \, \theta_1 -0.0581912473\,\sin \,\theta_1 )\\
 &+a_2 (0.1756420310\,\cos \, \theta_2 -0.3011001793\,\sin \, \theta_2 )\\
 &+a_3 ( 0.3250998384\,\cos \, \theta_3 +0.0519822050 \, \sin \, \theta_3)\\
 \Xi_{I} =& 0.1147989903,\\
 A_x= & 0,\\
 P_x(\theta)= &
   a_1 ( 0.0490168560\, \cos \, \theta_1 -0.8643742388\, \sin \, \theta_1 )\\
 &+a_2 ( 2.9257472656\, \cos \, \theta_2 + 0.5387775767 \, \sin \, \theta_2) \\
 &+a_3 ( 0.7848263649 \, \cos \, \theta_3 + 2.0413493194 \, \sin \, \theta_3)
\end{array}
\end{equation}

\begin{figure}[h]
\begin{center}
\begin{tabular}{lll}
(a) & (b) \\
\includegraphics[width=7cm]{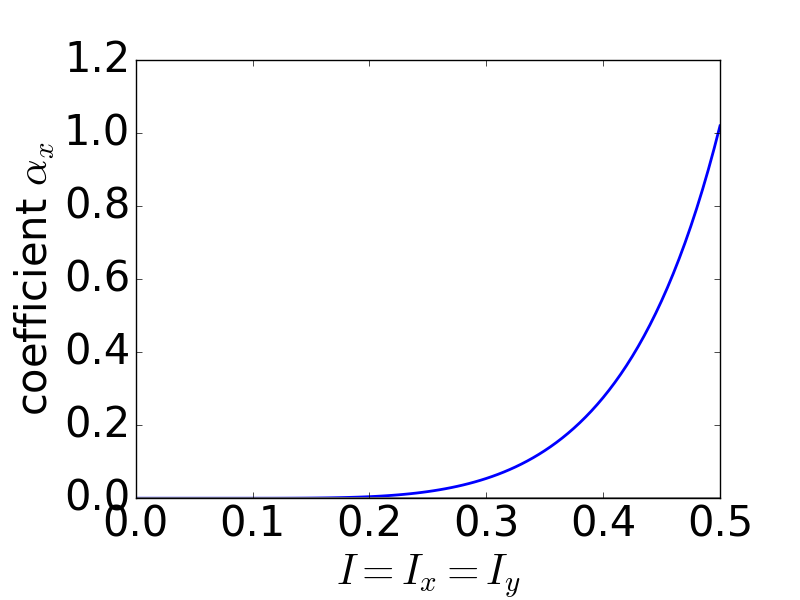} &
\includegraphics[width=7cm]{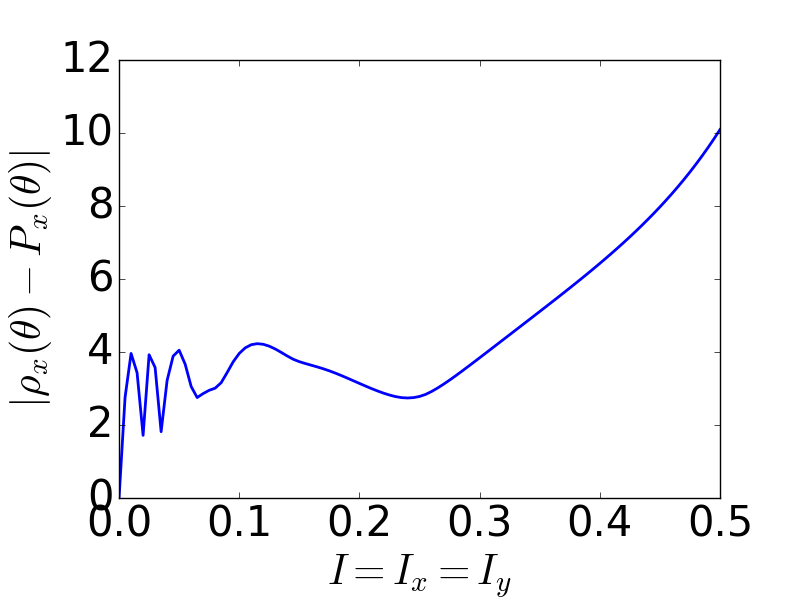} \\
\end{tabular}
\end{center}
\caption{
For different values of $I=I_x=I_y$ we show (a) the coefficient $\alpha_x$ in \eqref{eq:global_variacionals_ashwin_simplified} computed numerically using variational equations
along the heteroclinic connection to compare with the coefficient $A_x=0$ in \eqref{eq:AshwinMapMelnikovExample}
and (b) the difference in $L^1$-norm between the function $\rho_x(\theta)$ in \eqref{eq:global_variacionals_ashwin_simplified}
and the function $P_x(\theta)$ in \eqref{eq:AshwinMapMelnikovExample}.
}
\label{fig:comparison2_ashwin}
\end{figure}

We have checked that the global map
\eqref{eq:global_variacionals_ashwin_simplified} computed using variational equations, coincides with the map \eqref{eq:AshwinMapMelnikov} with the values given in \eqref{eq:AshwinMapMelnikovExample}
when $I_x$, $I_y$ are zero and it remains close as $I_x$ and $I_y$ increase.
To illustrate this, we show in Figure~\ref{fig:comparison2_ashwin} the
comparison between the coefficient $\alpha_x$ in \eqref{eq:global_variacionals_ashwin_simplified} and $A_x=0$ in
\eqref{eq:AshwinMapMelnikovExample}  and between the function $\rho_x(\theta)$ in \eqref{eq:global_variacionals_ashwin_simplified} and
$P_x(\theta)$ in \eqref{eq:AshwinMapMelnikovExample} as a function of $I=I_x=I_y$.

\section{Discussion}\label{sect:Discussion}

We have constructed the separatrix map for two different systems, the Duffing equation (see Section~\ref{sect:DuffingSeparatrix}) and
the 3-dimensional HBR model introduced in \cite{Ashwin2010}
(see Sections~\ref{sect:AshwinVar} and \ref{sect:AshwinMel}), both subject to periodic and, more relevantly for the
purpose of this study, quasi-periodic perturbations with at most 3 non-resonant frequencies.
The separatrix map associated to a Poincar\'{e} section constitutes a powerful tool to express in a simplified way the dynamics around homoclinic and heteroclinic trajectories of dynamical systems.
The ideas presented herein are extendable to more frequencies and larger networks.

We have obtained the technical results using two strategies: a Melnikov approach,
valid only for small perturbations but providing analytical descriptions,
and variational equations, based on analytico-numerical integration. We have compared both of them and checked its coincidence in the region where it is expected to be fulfilled.

For the Duffing equation, we have first analyzed the Lyapunov exponents for all the initial conditions. When perturbing with quasi-periodic perturbations with two or three non-resonant frequencies,
we generically obtain positive Lyapunov exponents, thus indicating chaotic behaviour.
We have then explored the distribution of the dominance times $\T_{dom}$ between impacts on Poincar\'{e} sections, whose histograms show a log-normal or Gamma distribution. The fitting quality obtained is comparable to that obtained from noisy perturbations of equivalent strength.

The results obtained for the Duffing equation persist for the HBR model. Remarkably, as for the Duffing equation,
from the separatrix map we also obtain a good agreement with Gamma distributions when perturbing with two and three non-resonant frequencies. For this case, we have developed a method based on
Melnikov integrals for non-Hamiltonian systems that have ``action-like'' variables vanishing at the heteroclinic connections.
This Melnikov approach, as opposed to the case for the Duffing equation, requires the numerical computation of the integrals because of the cumbersome parameterization of the non-perturbed heteroclinic orbits.
The resulting discrete model can be thought of as a new model for bistable perception, much easier to use than the full model expressed in terms of differential equations.
Altogether, the methodology that we propose provides in both cases an alternative discrete model (a map) which avoids numerically unstable computations. More precisely, time-continuous models require the numerical integration close to saddle points while the separatrix maps resolves this issue by using the linear approximation around the saddles. More refined maps could be obtained by substituting the local maps by higher order approximations (normal forms).

From a modelling point of view, we have proved that important features attributed to psychophysical experiments of bistable perception, namely the Gamma distribution of dominance times,
cannot only be reproduced by noisy perturbations but also by quasi-periodic perturbations with two or more non-resonant frequencies.
This fact was known for noisy perturbations but not for deterministic ones
\cite{Moreno2007,Ashwin2010,Pastukhov2013}.
It is worth noticing that the signal of a noisy perturbation presents a continuous spectrum and so, our result implies that the same output distribution can be achieved by perturbing the system with only few frequencies. One could argue that for finite time simulations, as for instance those in \cite{Ashwin2010}, the spectrum is less richer than the theoretical prediction for infinite time,
but still the support of the spectrum in the input distributions decays drastically in size when jumping from noisy to quasi-periodic perturbations.

We would also like to draw the attention to the question of what the noise is actually
representing, since models in the literature are
not precise enough about the source of the stochastic nature of bistable perception. It is believed that perceptual switches are spontaneous and stochastic events (for instance, {\it a priorities}) which cannot be eliminated by intentional efforts and it has been largely emphasized the relative role of noise versus adaptation \cite{Moreno2007,Kang2010} but, as far as we know, there are no solid arguments that sustain that they must be forcedly spontaneous and purely stochastic. The models usually contain two main variables that represent the ensembles of neurons more directly related to the percepts (for instance, the left and right eyes in binocular rivalry, the most well-known phenomenon of bistable perception). Adding noise to the model, on the other hand, entails the assumption that this basic biperceptual system is receiving inputs from a continuous spectrum. However,
Electroencephalography (EEG) studies (see for instance \cite{Doesburg2005,Doesburg2009}) suggest a prevalent role of gamma-band frequencies.
More precisely, transient gamma-band synchrony in localized recurrent prefrontal and parietal brain areas (responsible for executive functions) have been reported to precede switching between percepts in binocular rivalry. These findings make plausible the conjecture that a few number of frequencies in the input sources could be sufficient to account for the statistics of perceptual switches.

This perspective of bistable perception using maps brings up new possibilities to investigate this phenomenon which are beyond the scope of this paper. For instance, studying in depth the dynamics of these maps or its fitting to experimental data. Here, we give a first step in this direction by computing the Lyapunov exponents and certifying the compatibility with the obtained Gamma distributions. As a future work, we plan to use the separatrix map models to fit other experimental (psychophysical) data.

We finish by pointing out that our results extend naturally to other problems modelled by means of heteroclinic networks already mentioned in the Introduction like decision-making, memory-retrieval, central patterns generators or ecological models.

\section{Acknowledgments}
This work has been partially funded by the Spanish grants MINECO-FEDER-UE MTM2015-65715-P (AD, GH), MINECO-FEDER-UE MTM-2015-71509-C2-2-R (AG,GH), the Catalan Grant 2017SGR1049 (AD, AG, GH).
GH wants to acknowledge the RyC project RYC-2014-15866. We thank A. Viero for providing us information about MEGNO.
We also acknowledge the use of the UPC Dynamical Systems group's
cluster for research computing (https://dynamicalsystems.upc.edu/en/computing/).

\bibliographystyle{abbrv}

\bibliography{references}

\appendix

\section*{Appendices}

\section{Variational equations}\label{ap:variational}
Consider the system of the first variational equations along the separatrix $\Gamma$ for the extended
system, consisting of system~\eqref{eq:systemuv} with the extra equation $\dot{\eps}=0$, as introduced in Section~\ref{sect:VariationalDuffing},
given by
\[\frac{d}{dt} D_\wv \hat \varphi (t;\wv)=A(t) D_{\wv} \hat \varphi (t;\wv), \qquad D_{\wv} \hat \varphi (0;\wv)=Id_{n+3},\]
where
\[
A(t)=
 \left (
 \begin{array}{cccc}
 \dfrac{\partial F^u}{\partial u} & \dfrac{\partial F^u}{\partial v} & 0 & \dfrac{\partial F^u}{\partial \eps} \vspace{0.2cm} \\
 \dfrac{\partial F^v}{\partial u} & \dfrac{\partial F^v}{\partial v} & 0 & \dfrac{\partial F^v}{\partial \eps} \vspace{0.2cm}\\
 0 & 0 & 0 & 0 \\
 0 & 0 & 0 & 0 \\
 \end{array}
 \right )_{\mbox{\normalsize $| \, \hat \varphi (t; \wv^s)$}}
\]
with $F^u,F^v$ given in \eqref{eq:systemuv}, $\hat \varphi$ is the flow of the extended system and
$\wv^s=(u^s,v^{*},\theta^s,0)$ with $(u^s,v^*)=\{ \Gamma \cap J^{out} \}$.

Let us denote $\hat \varphi_w$
the derivative with respect to $w$ and $\hat \varphi^w$ the coordinate $w$, for $w=u,v,\theta,\eps$.
Since
\[
\fv_{u,v,\theta,\eps}^{\theta}=0, \qquad \textrm{then} \qquad \hat \varphi^{\theta}_{u,v,\eps}=0, \, \hat \varphi^{\theta}_{\theta}=1, \\
\]
and since
\[
\fv_{u,v,\theta,\eps}^{\eps}=0, \qquad \textrm{then} \qquad \hat \varphi^{\eps}_{u,v,\theta}=0, \, \hat \varphi^{\eps}_{\eps}=1. \\
\]
Therefore, we are left with the following equations
\[
\begin{array}{rclcrcl}
\fv_u^u & = &  \dfrac{\partial F^u}{\partial u} \hat \varphi_u^u + \dfrac{\partial F^u}{\partial v} \hat \varphi_u^v, & \qquad&
\fv_u^v & = &  \dfrac{\partial F^v}{\partial u} \hat \varphi_u^u + \dfrac{\partial F^v}{\partial v} \hat \varphi_u^v, \vspace{0.2cm}\\
\fv_v^u & = &  \dfrac{\partial F^u}{\partial u} \hat \varphi_v^u + \dfrac{\partial F^u}{\partial v} \hat \varphi_v^v, & \qquad &
\fv_v^u & = &  \dfrac{\partial F^u}{\partial u} \hat \varphi_v^u + \dfrac{\partial F^u}{\partial v} \hat \varphi_v^v, \\
\fv_{\theta}^u & = &  0, & \qquad &
\fv_{\theta}^v & = &  0, \\
\fv_{\eps}^u & = &  \dfrac{\partial F^u}{\partial u} \hat \varphi_{\eps}^u + \dfrac{\partial F^u}{\partial v} \hat \varphi_{\eps}^v + \dfrac{\partial F^u}{\partial \eps}, & \qquad &
\fv_{\eps}^u & = &  \dfrac{\partial F^u}{\partial u} \hat \varphi_{\eps}^u + \dfrac{\partial F^u}{\partial v} \hat \varphi_{\eps}^v + \dfrac{\partial F^u}{\partial \eps}. \\
\end{array}
\]

Notice that only the equations for $\hat \varphi^u_{\eps}$ and $\hat \varphi^v_{\eps}$ depend on $\theta$
through the term $\partial F^{u}/\partial \eps$ and $\partial F^{v}/\partial \eps$, respectively. Therefore,
we will compute $\hat \varphi^u_{\eps}$ and $\hat \varphi^v_{\eps}$ for different initial conditions of
$\theta^s$. We use Fourier series to obtain an analytical expression for $\hat \varphi^u_{\eps}$ and $\hat \varphi^v_{\eps}$
as a function of $\theta$.

From the solution of the variational equations, we obtain the first order approximation of the global map $\hat T_G$
(see equation~\eqref{eq:Tgbar}).
Indeed, let us take a vector
$(\Delta u,0,\Delta \theta, \eps)$ onto the Poincar\'{e} section $\Sigma^{out}$ (where $v=v^*$, see
\eqref{eq:section_def}) and compute $\alpha_v, \alpha_{\theta},\kappa_v,\kappa_{\theta},\rhoF_v,\rhoF_{\theta}$ such that:
\[
\left (
\begin{array}{cccc}
\hat \varphi_u^u & \hat \varphi_v^u & 0 & \hat \varphi_{\eps}^u \\
\hat \varphi_u^v & \hat \varphi_v^v & 0 & \hat \varphi_{\eps}^v \\
0 & 0 & 1 & 0 \\
0 & 0 & 0 & 1
\end{array}
\right )
_{|(\tau^*;\wv^s)}
\left (
\begin{array}{c}
\Delta u \\
0 \\
\Delta \theta \\
\eps
\end{array}
\right )
+
\Delta t
\left (
\begin{array}{c}
F^u \\
F^v \\
\omega \\
0
\end{array}
\right )
_{|\hat \varphi (\tau^*;\wv^s)}
=
\left (
\begin{array}{c}
0 \\
\alpha_v \Delta u + \kappa_v \Delta \theta + \rhoF_v (\theta) \eps \\
\alpha_{\theta} \Delta u + \kappa_{\theta} \Delta \theta + \rhoF_{\theta} (\theta) \eps \\
\eps
\end{array}
\right ),
\]
where $\tau^*=\tau^*(\wv^s)$ is such that $\hat \varphi^u (\tau^*(\wv^s);\wv^s)=u^*$.
From the first coordinate we obtain
\[ \Delta t = - \frac{\hat \varphi_u^u \Delta u + \hat \varphi_{\eps}^u \Delta \eps}{F^u},\]
and therefore,
\begin{equation}\label{eq:constants_ap}
\begin{array}{ccc}
\alpha_v = \hat \varphi_u^v- \dfrac{F^v}{F^u} \hat \varphi_u^u, & \kappa_v= 0, & \rhoF_v (\theta)= \hat \varphi_{\eps}^v- \dfrac{F^v}{F^u} \hat \varphi_{\eps}^u, \vspace{0.2cm}\\
\alpha_{\theta} = - \omega \dfrac{\hat \varphi_{u}^u}{F^u}, & \kappa_{\theta} = 1, & \rhoF_{\theta}(\theta)= - \omega\dfrac{\hat \varphi_{\eps}^u}{F^u}, \\
\end{array}
\end{equation}
obtaining the formulas that are given in \eqref{eq:Tgbar} and \eqref{eq:coef_uv}.

\section{Choice of the Poincar\'{e} sections}\label{ap:cr}

\begin{figure}[h]
\begin{center}
\begin{tabular}{ll}
(a) Comparison scheme & (b) $\gamma=0.008$ \\
\includegraphics[width=6cm]{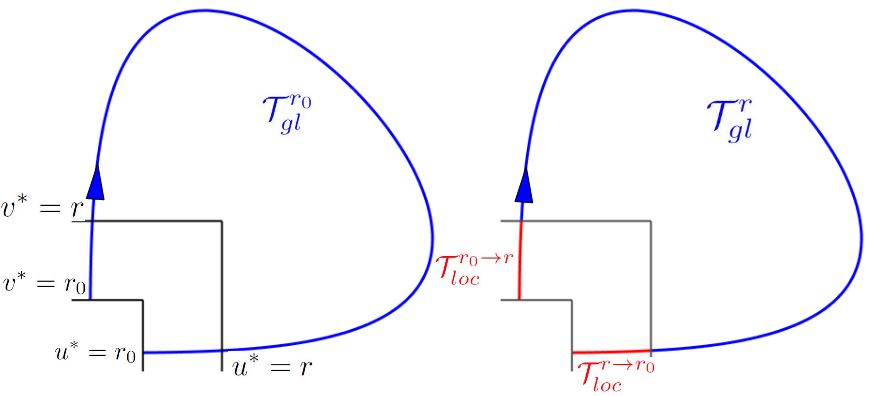} &
\includegraphics[width=6cm]{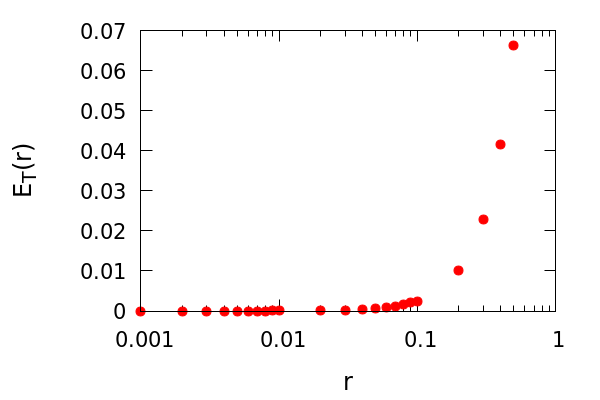} \\
(c) $\gamma=0.08$ & (d) HBR model \\
\includegraphics[width=6cm]{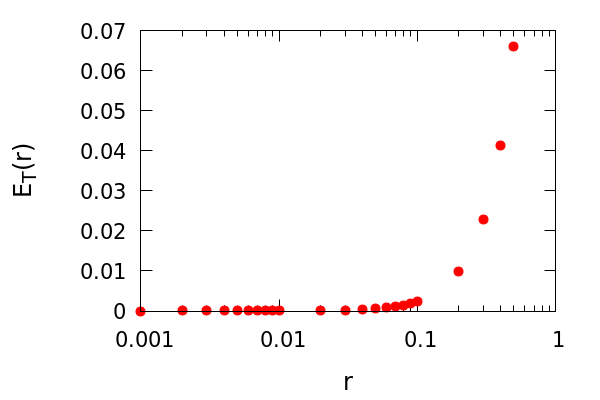} &
\includegraphics[width=6cm]{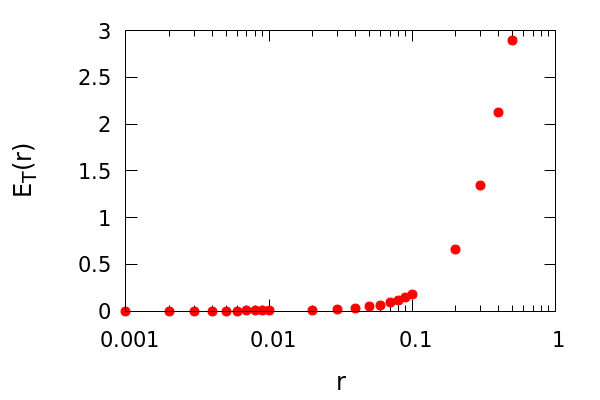} \\
\end{tabular}
\end{center}
\caption{Justification for the choice of $r$. (a) Schematic representation of the two compared times: global map (left) and concatenation of three maps (right). On trajectories shown in blue, the time is computed using global maps, while on trajectories shown in red, the time is computed using local approximations. (b,c,d) Function $E_{\T}(r)$ showing the difference between the time from section
$\Sigma^{out}_{0.001}$ to $\Sigma^{in}_{0.001}$ computed using a global map or a combination of
local and global maps involving $r$  for (a) the Duffing equation with $\gamma=0.008$,
(b) the Duffing equation with $\gamma=0.08$ and (c) the HBR model
(see equation \eqref{eq:errort} and Appendix~\ref{ap:cr} for more details).
}
\label{fig:comparison_errors}
\end{figure}

In both examples, the Duffing equation and the HBR model, we have chosen sections $\Sigma^{in}$ and
$\Sigma^{out}$ located at a distance $r=0.1$ from the saddle points. This choice is a compromise between avoiding long computations along the separatrices, which requires $r$ to be as large as possible, and maintaining the validity of the approximation of the local dynamics by the linear map, which requires $r$ to be as small as possible. To assess this balance, we considered fixed sections at a very small distance ($r_0$) to the saddle and, for $r>r_0$, we evaluated the time error induced by the fact of considering the approximation of the local dynamics from section $r_0$ to section $r$ instead of considering the global map.
More precisely, for the Duffing equation (for the HBR model, the procedure is similar), we compare the global time from section $v=r_0$ to $u=r_0$ (denoted by $\T_{gl}^{r_0}$) along the separatrix with the concatenation of three times (see also Figure~\ref{fig:comparison_errors}(a)):
(1) the time $\T_{loc}^{r_0 \rightarrow r}$ to go from section $v=r_0$ to
section $v=r>r_0$ computed using the local approximation ($\T_{loc}^{r_0 \rightarrow r}=1/\lambda_{+} \ln(r/r_0)$)
plus (2) the global time ($\T_{gl}^{r}$) along the separatrix to go from section $v=r$ to
$u=r$ plus (3) the time $\T_{loc}^{r \rightarrow r_0}$ to go from section
$u=r>r_0$ to section $u=r_0$ computed using the local approximation ($\T_{loc}^{r \rightarrow r_0}=1/\lambda_{-} \ln(r_0/r)$). That is,
we compute the error function
\begin{equation}\label{eq:errort}
E_\T(r)= |\T_{gl}^{r_0} - (\T_{loc}^{r_0 \rightarrow r} + \T_{gl}^{r} +\T_{loc}^{r \rightarrow r_0})|.
\end{equation}
In Figures~\ref{fig:comparison_errors}(b) and (c) we show the function
$E_\T(r)$ with $r_0=0.001$ and $0.001 \leq r \leq 0.5$ for $\gamma=0.008$ and $\gamma=0.08$, respectively.
By looking at these plots, it is clear that $r=0.1$ is a good compromise.

An analogous computation for the HBR model gives the results shown
in Figure~\ref{fig:comparison_errors}(d). Again it is clear that $r=0.1$ is a good compromise also for this model.

\section{Fitting of the histograms}\label{ap:fitting}

The distributions of dominance times in Figures~\ref{fig:hist_gae-3}, \ref{fig:hist_gae-2}, \ref{fig:hist_noise_uv} (top)
and \ref{fig:DominanceAshwin} top  have been fitted to log-normal and Gamma distributions.
The probability density function for the log-normal distribution is
\[
f_{ln}(x)=\frac{1}{\sigma x \sqrt {2 \pi}} \exp \left (- \frac{(\ln(x)-\mu)^2}{2\sigma^2} \right ),
\]
where $\mu$ is the mean and $\sigma$ is the standard deviation of the normally distributed logarithm of the variable.
The probability density function for the Gamma distribution with a shape parameter $a$
and a scale parameter $\lambda$ is
\[
f_{g}(x) = \frac{1}{\Gamma(a)  \lambda^a } x^{a-1} \exp(-x/\lambda).
\]

Maximum likelihood fits of the time distributions to a log-normal and Gamma distribution give the
parameter values indicated in Table~\ref{tab:dt_duffing} for the Duffing equation and Table~\ref{tab:dt_ashwin}
for the HBR model.

The distributions of impacts on sections $\Sigma^{out}$ in Figures~\ref{fig:hist_gae-3}, \ref{fig:hist_gae-2}, \ref{fig:hist_noise_uv} (bottom)
have been fitted to a normal distribution. The probability density function for the normal distribution is
\[
f_{n}(x)=\frac{1}{\sigma \sqrt {2 \pi}} \exp \left (- \frac{(x-\mu)^2}{2\sigma^2} \right ),
\]
where $\mu$ is the mean and $\sigma$ is the standard deviation. The parameter values obtained from the maximum likelihood fits to the normal distribution
for the Duffing equation are given in Table~\ref{tab:dt_duffing}.

\begin{table}
\begin{center}
\renewcommand{\arraystretch}{1.5}
\begin{tabular}{cc|cccc|cc}
\multicolumn{2}{c}{} & \multicolumn{2}{c}{Gamma} & \multicolumn{2}{c}{Log-normal} & \multicolumn{2}{c}{Normal}\\
\cline{2-8}
& Model       &  shape ($a$) & scale ($\lambda$) & $\sigma$ & $\mu$ (scale=$e^\mu$) & $\mu$ & $\sigma$\\
\cline{2-8}
\multirow{3}{*}{\begin{sideways}$\gamma=0.008$\end{sideways}}
& SM 2 freq & 105.584 & 0.08686 & 0.09766 & 0.96084 & 0.00022 & 0.02664 \\
& SM 3 freq & 87.2084 & 0.09990 & 0.10748 & 0.93821 & -0.00030 & 0.04286 \\
&noise & 84.2724 & 0.10647 & 0.10938 & 0.95094 & -0.00073 & 0.03466 \\
\cline{2-8}
\multirow{3}{*}{\begin{sideways}$\gamma=0.08$\end{sideways}} &
 SM 2 freq & 126.606 & 0.07405 & 0.08922 & 0.97069 & 0.00026 & 0.01996  \\
& SM 3 freq & 98.6673 & 0.09110 & 0.10105 & 0.95199 & 0.00026 & 0.03072 \\
&noise & 95.3895 & 0.09917 & 0.10273 & 0.97411 & -0.00093 & 0.01932 \\
\cline{2-8}
\end{tabular}
\caption{Parameters of the Gamma and log-normal distributions obtained from the fitting of the dominance times histograms and
parameters of the normal distribution obtained from the fitting of the impact histograms shown in Figures~\ref{fig:hist_gae-3}, \ref{fig:hist_gae-2}, \ref{fig:hist_noise_uv}
for simulations of the separatrix map (SM) with 2 and 3 frequencies and the system with noise corresponding to the Duffing equation.}
\label{tab:dt_duffing}
\end{center}
\end{table}

\begin{table}
\begin{center}
\renewcommand{\arraystretch}{1.5}
\begin{tabular}{c|cccc}
\multicolumn{1}{c}{} & \multicolumn{2}{c}{Gamma} & \multicolumn{2}{c}{Log-normal} \\
\hline
Model &  shape ($a$) & scale ($\lambda$) & $\sigma$ & $\mu$ (scale=$e^\mu$)  \\
\hline
SM 1 freq & 167.492 & 0.42598  & 0.07735 & 4.26527 \\
SM 2 freq & 87.1589 & 0.65399 & 0.10765 & 4.03872 \\
SM 3 freq & 55.9174 & 1.00719  & 0.13446 & 4.02414 \\
noise & 38.0614 & 1.49805 & 0.16332 & 4.03324 \\
\end{tabular}
\caption{Parameters of the Gamma and log-normal distributions obtained from the fitting of the dominance times histograms shown in
Figure~\ref{fig:DominanceAshwin} for simulations of the separatrix map (SM) with 1, 2 and 3 frequencies and the system with noise corresponding to the HBR model.}
\label{tab:dt_ashwin}
\end{center}
\end{table}

\pagebreak

\section*{List of Figures}


\end{document}